\documentclass[aop,MSNbibl,dvips]{arximspdf}
\usepackage{graphicx,breakurl}
\doi{10.1214/13-AOP864} 
\volume{42}
\issue{4}
\pubyear{2014}
\firstpage{1396}
\lastpage{1437}
\makeatletter
\newproclaim{step}{Step}
\newproclaim{stepo}{Step}
\newcommand{\eps}{\varepsilon}
\newcommand{\Ii}{\mathcal{I}}
\newcommand{\cov}{\mathop{\operatorname{Cov}}}
\newcommand{\orb}{\mathop{\operatorname{Orb}}}
\newcommand{\E}{\mathbf{E}}
\renewcommand{\P}{\mathbf{P}}
\newcommand{\NN}{\mathbb{N}}
\newcommand{\tr}{\mathop{\operatorname{tr}}}
\newcommand{\eqd}{\sim}
\newcommand{\Exp}{\operatorname{Exp}}
\newcommand{\leb}{\operatorname{Leb}}
\newcommand{\limitN}{{N}}
\newcommand{\Ww}{\mathcal{W}}
\newcommand{\Poi}{\operatorname{Poi}}
\def\mid{|}
\newcommand{\divides}{\mid}
\newcommand{\Tt}{\mathcal{T}}
\newcommand{\Y}{\mathbf{Y}}
\newcommand{\Aa}{\mathcal{A}}
\newcommand{\ZZ}{\mathbb{Z}}
\newcommand{\RR}{\mathbb{R}}
\renewcommand{\d}{\mathrm{d}}
\newcommand{\eqref}[1]{(\ref{#1})}
\newtheorem{thmm}{Theorem}
\newtheorem{lemma}[thmm]{Lemma}
\newtheorem{prop}[thmm]{Proposition}
\newtheorem{cor}[thmm]{Corollary}
\newproclaim{rmk}[thmm]{Remark}
\newproclaim{exm}[thmm]{Example}
\newproclaim{clm}[thmm]{Claim}
\newproclaim{defn}[thmm]{Definition}
\makeatother

\begin{document}
\begin{frontmatter}

\title{Cycles and eigenvalues of sequentially growing random regular
graphs\thanksref{T1}}
\runtitle{Cycles and eigenvalues}
\thankstext{T1}{Supported in part by NSF Grant DMS-10-07563.}

\begin{aug}
\author[A]{\fnms{Tobias} \snm{Johnson}\corref{}\ead[label=e1]{toby@math.washington.edu}}
\and
\author[A]{\fnms{Soumik} \snm{Pal}\ead[label=e2]{soumik@u.washington.edu}}
\runauthor{T.~Johnson and S.~Pal}
\affiliation{University of Washington}
\address[A]{Department of Mathematics\\
University of Washington\\
Seattle, Washington 98195\\
USA\\
\printead{e1}\\
\phantom{E-mail:\ }\printead*{e2}} 
\end{aug}

\received{\smonth{5} \syear{2012}}
\revised{\smonth{5} \syear{2013}}

%
\begin{abstract}
Consider the sum of $d$ many i.i.d. random permutation matrices on $n$
labels along with their transposes. The resulting matrix is the
adjacency matrix of a random regular (multi)-graph of degree $2d$ on
$n$ vertices. It is known that the distribution of smooth linear
eigenvalue statistics of this matrix is given asymptotically
by sums of Poisson random variables. This is in contrast with Gaussian
fluctuation of similar quantities in the case of Wigner matrices. It is
also known that for Wigner matrices the joint fluctuation of linear
eigenvalue statistics across minors of growing sizes can be expressed
in terms of the Gaussian Free Field (GFF). In this article, we explore
joint asymptotic (in~$n$) fluctuation for a coupling of all random
regular graphs of various degrees obtained by growing each component
permutation according to the Chinese Restaurant Process. Our primary
result is that the corresponding eigenvalue statistics can be expressed
in terms of a family of independent Yule processes with immigration.
These processes track the evolution of \mbox{short} cycles in the graph.
If we now take $d$ to infinity, certain GFF-like properties emerge.
\end{abstract}

%
\begin{keyword}[class=AMS]
\kwd{60B20}
\kwd{05C80}
\end{keyword}
\begin{keyword}
\kwd{Random regular graphs}
\kwd{eigenvalue fluctuations}
\kwd{Chinese restaurant process}
\kwd{minors of random matrices}
\end{keyword}

\end{frontmatter}

\section{Introduction} We consider graphs that have labeled vertices
and are regular, that is, every vertex has the same degree.
We allow our graphs to have loops and multiple edges (such graphs
are sometimes called multigraphs or pseudographs). 
Additionally, our graphs will be sparse in the sense that the degree
will be negligible compared to the order.
Every such graph has an associated adjacency matrix whose $(i,j)$th
element is the number of edges between vertices $i$ and $j$, with loops
counted twice.
When the graph is randomly selected, the matrix is random, and we are
interested in studying the eigenvalues of the resulting symmetric
matrix. Note that, due to regularity, it does not matter whether we
consider the eigenvalues of the adjacency or the Laplacian matrix.

The precise distribution of this random regular graph is somewhat ad
hoc. We will use what is called the permutation model. Consider the
permutation digraphs generated by $d$ many i.i.d. random permutations
on $n$ labels. We remove the direction of the edge and \textit
{collapse} all these graphs on one another. This results in 
a $2d$-regular graph on $n$ vertices, denoted by $G(n,2d)$. At the
matrix level this is given by adding all the $d$ permutation matrices
and their transposes.

Our present work is an extension of the study of eigenvalue
fluctuations carried out in \cite{DJPP}. We are
motivated by the recent work by Borodin on joint eigenvalue
fluctuations of minors of Wigner matrices and the (massless or
zero-boundary) Gaussian Free Field (GFF) \cite{Bor1,Bor2}. Eigenvalues
of minors are closely related to interacting particle systems \cite
{Fer10,FF10}, and the KPZ universality class of random surfaces \cite
{BorFer}. See \cite{JN} for more on eigenvalues of minors of GUE and
\cite{ANvM} for those of Dyson's Brownian motion.

Let us consider a particular but important case of Borodin's result in
\cite{Bor1} (single sequence, the entire $\NN$). An $n \times n$ real
symmetric Wigner matrix has i.i.d. upper triangular off-diagonal
elements with four moments identical to the standard Gaussian. The
diagonal elements are usually taken to be i.i.d. with mean zero
variance two. Notice that every principal submatrix (called minors in
this context) of a Wigner matrix is again a Wigner matrix of a smaller
order. Thus, on some probability space one can construct an infinite
order Wigner matrix $W$ whose $n\times n$ minor $W(n)$ is a Wigner
matrix of order $n$.

Let $z$ be a complex number in the upper half plane $\mathbb{H}$.
Define $y=\vert z\vert^2$ and $x=2 \Re(z)$.
Consider the minor $W(\lfloor ny \rfloor)$, and
let $N(z)$ be the number of its eigenvalues that are
greater than or equal to $\sqrt{n}x$.
Define the \textit{height function}
%
%
\begin{equation}
\label{eq:heightfn} H_n(z):= \sqrt{\frac{\pi}{2}} N(z).
\end{equation}
Then Borodin shows that $\{ H_n(z) - E H_n(z), z\in\mathbb{H} \}$,
viewed as distributions, converges in law to a generalized Gaussian
process on $\mathbb{H}$ with a covariance kernel
%
%
\begin{equation}
\label{eq:GFFcov} C(z, w)= \frac{1}{2\pi} \ln\biggl\vert\frac{z- w}{z - \overline{w}} \biggr\vert.
\end{equation}
The above is the covariance kernel for the GFF on the upper half plane.

An equivalent assertion is the following. Let $[n]$ denote the set of
integers $\{1,2,\ldots, n\}$. Consider the Chebyshev polynomials of the
first kind,
$\{ T_n, n=0,1,2,\ldots\}$, on the interval $[-1,1]$. These
polynomials are given by the identity $T_n(\cos(\theta))\equiv\cos
(n\theta)$. We specialize \cite{Bor1}, Proposition 3, for the case of
GOE ($\beta=1$).
Fix $m$ positive real numbers $t_1 < t_2 < \cdots< t_m$. In the
notation of \cite{Bor1}, we take $L=n$ and $B_i(n)=[\lfloor t_i n
\rfloor]$. Then, for any positive integers $j_1, j_2, \ldots, j_m$, the
random vector
\[
\bigl( \tr T_{j_i} \bigl( W \bigl(\lfloor t_i n \rfloor
\bigr)/2\sqrt{t_i n} \bigr) - \E\tr T_{j_i} \bigl( W \bigl(
\lfloor t_i n \rfloor \bigr) /2\sqrt{t_i n} \bigr), i
\in[m] \bigr)
\]
converges in law, as $n$ tends to infinity, to a centered Gaussian vector.
For $s\leq t$,
%
%
\begin{equation}
\label{eq:chebycovBor}\quad  \lim_{n\rightarrow\infty} \cov\bigl( \tr
T_i \bigl( W\bigl( \lfloor t n \rfloor\bigr)/ 2\sqrt{tn} \bigr), \tr
T_k \bigl( W \bigl(\lfloor s n \rfloor \bigr) /2\sqrt{sn} \bigr)
\bigr)=\delta_{ik} \frac{k}{2} \biggl(\frac{s}{t}
\biggr)^{k/2},
\end{equation}
which gives the covariance kernel of the limiting vector.
In particular, all such covariances are zero when $i\neq k$. Note that
the traces can be expressed as integrals of the height function of the
corresponding submatrices. Thus, by approximating continuous compactly
supported functions of $z$ by a function that is piecewise constant in
$y$ and polynomial in $x$, one gets the kernel \eqref{eq:GFFcov}.

\subsection{Main results}\label{sec:mainresults} By a tower of random
permutations, we mean a sequence of random permutations $(\pi^{(n)},
n\in\NN)$ such that:
\begin{longlist}[(ii)]
\item[(i)] $\pi^{(n)}$ is a uniformly distributed random permutation of
$[n]$ for each $n$, and
\item[(ii)] for each $n$, if $\pi^{(n)}$ is written as a product of
cycles then $\pi^{(n-1)}$ is derived from $\pi^{(n)}$ by deletion of
the element $n$ from its cycle.
\end{longlist}
The stochastic process that grows $\pi^{(n)}$ from $\pi^{(n-1)}$ by
sequentially inserting an element $n$ randomly is called the Chinese
Restaurant Process (CRP).
We will review the basic principles at a later section. In \cite{KOV}
and other related work,
a sequence of permutations satisfying condition (ii) is called
a \emph{virtual permutation}, and the distribution on virtual permutations
satisfying condition (i) is considered as a substitute for Haar
measure on $S(\infty)$,
the infinite symmetric group. This is used to study the representation
theory of $S(\infty)$, with connections to random matrix theory.
A~recent extension of this idea is
\cite{BNN}.

Now suppose we construct a countable collection $\{ \Pi_d, d \in\NN
\}$ of towers of random permutations. We will denote the permutations
in $\Pi_d$ by $\{ \pi_d^{(n)}, n \in\NN\}$. Then it is possible to
model every possible $G(n, 2d)$ by adding the permutation matrices (and
their transposes) corresponding to $\{ \pi^{(n)}_j, 1\le j \le d \}$.
In what follows, we will keep $d$ fixed and consider $n$ as a growing
parameter. Thus, $G_n$ will represent $G(n,2d)$ for some fixed $d$.
Here and later, $G_0$ will represent the empty graph.
We construct a continuous-time version of this by inserting new vertices
into $G_n$ with rate $n+1$. Formally,
define independent \label{page:poissonization}
times $T_i\eqd\Exp(i)$, and let
\[
M_t=\max \Biggl\{m\dvtx \sum_{i=1}^mT_i
\leq t \Biggr\},
\]
and define the continuous-time Markov chain
$G(t)=G_{M_t}$.\label{page:chaindef}
When $d=1$, this process is essentially just a continuous-time version
of the CRP itself. Though this case is unusual compared to the rest---for
example, $G(t)$ is likely to be disconnected when $d=1$ and
connected when $d$ is larger---our results do still hold.

Our first result is about the process of \textit{short cycles} in the
graph process $G(t)$.
By a cycle of length $k$ in a graph, we mean what is sometimes called
a simple cycle: a walk in the graph that begins
and ends at the same vertex, and that otherwise repeats no vertices.
We will give a more formal definition in Section~\ref{sec:wordcombinatorics}.
Let $(C_k^{(s)}(t), k \in\NN)$ denote the number of cycles of
various lengths $k$ that are present in $G(s+t)$. This process is not
Markov, but nonetheless it converges to a Markov process (indexed by
$t$) as $s$ tends to infinity. 

To describe the limit, define 
\[
a(d,k)=\cases{(2d-1)^{k} - 1 + 2d,& \quad $\mbox{when $k$ is even}$,
\vspace*{2pt}
\cr
(2d-1)^{k} +1,&\quad  $\mbox{when $k$ is odd}.$}
\]
Consider the set of natural numbers $\NN=\{1,2,\ldots\}$ with the measure
\[
\mu(k)= \tfrac{1}{2} \bigl[ a(d,k) - a(d,k-1) \bigr], \qquad k\in\NN, a(d,0):=0.
\]
%
Consider a Poisson point process $\chi$ on $\NN\times[0, \infty)$
with an intensity measure given on $\NN\times(0, \infty)$
by the product measure $\mu\otimes\leb$, where $\leb$ is the Lebesgue
measure, and with additional masses of $a(d,k)/2k$ on $ (k,0)$ for
$k\in\NN$.

Let $\widetilde{P}_{x}$ denote the law of an one-dimensional pure-birth
process on $\NN$ given by the generator:
\[
L f(k) = k \bigl( f(k+1) - f(k) \bigr), \qquad k \in\NN,
\]
starting from $x\in\NN$. This is also known as the \textit{Yule process}.

Suppose we are given a realization of $\chi$. For any atom $(k,y)$ of
the countably many atoms of $\chi$, we start an independent process
$(X_{k,y}(t), t\ge0)$ with law $\widetilde{P}_{k}$. Define the random sequence
\[
\limitN_k(t):= \sum_{(j,y)\in\chi\cap\{ [k] \times[0, t] \} } 1 \bigl\{
X_{j,y}(t-y)=k \bigr\}.
\]
In other words, at time $t$, for every site $k$,
we count how many of the processes that started at time $y\le t$ at
site $j \le k$ are currently at $k$. 
Note that both $(\limitN_k(\cdot), k\in\NN)$ and $(\limitN
_k(\cdot
), k\in[K])$, for some $K\in\NN$, are Markov processes, while
$N_k(\cdot)$ for
fixed $k$ is not.

%
%
\begin{thmm}\label{mainthm:cycles}
As $s\to\infty$, the process $(C_k^{(s)}(t), k\in\NN, 0\le t <
\infty)$ converges in law
in $D_{\mathbb{R}^{\infty}}[0, \infty)$ to the Markov
process $(\limitN_k(t), k\in\NN, 0\le t < \infty)$. The limiting
process is stationary.
\end{thmm}

%
\begin{rmk}
In fact, the same argument used to prove Theorem \ref{mainthm:cycles}
shows that the process $(C_k^{(s)}(t), -\infty< t < \infty)$
converges in law to the Markov process $(\limitN_k(t), -\infty< t <
\infty)$ running in stationarity. The same conclusion holds for all the
following theorems in this section.
\end{rmk}

We now explore the joint convergence across various $d$'s. Define
$C_{d,k}^{(s)}(t)$ naturally, stressing the dependence on the parameter $d$.

%
\begin{thmm}\label{mainthm:intertwining}
There is a joint process convergence of $(C_{i,k}^{(s)}(t), k \in
\NN
, i \in[d], t \ge0)$ to a limiting process $(\limitN_{i,k}(t),
k \in\NN, i\in[d], t\ge0)$. This limit is a Markov process
whose marginal law for every fixed $d$ is described in Theorem~\ref
{mainthm:cycles}.
Moreover, for any $d\in\NN$, the process $(\limitN_{d+1, k}(\cdot) -
\limitN_{d,k}(\cdot), k \in\NN)$ is independent of the process
$(\limitN_{i,k}(\cdot), k \in\NN, i\in[d])$ and evolves as a
Markov process.
Its generator (defined on functions dependent on finitely many
coordinates) is given by
\[
Lf(x)= \sum_{k=1}^\infty k x_k
\bigl[ f ( x + e_{k+1} - e_k ) - f(x) \bigr] + \sum
_{k=1}^\infty\nu(d,k) \bigl[ f(x + e_k) -
f(x) \bigr],
\]
where $x$ is a nonnegative sequence, $(e_k, k \in\NN)$ are the
canonical orthonormal basis of $\ell^2$, and
\[
\nu(d,k)=\tfrac{1}{2} \bigl[ a(d+1, k) - a(d+1,k-1) - a(d,k) + a(d, k-1)
\bigr].
\]
\end{thmm}

%
\begin{rmk}
Theorems \ref{mainthm:cycles} and \ref{mainthm:intertwining} show an
underlying branching process structure. We actually prove a more
general decomposition where cycles are tracked by edge labels. The
additive structure also imparts a natural intertwining relationship
between the Markov operators. See \cite{CPY}, Section~2 and \cite
{DF90,Bor1}.
\end{rmk}

We now focus on eigenvalues of $G(t)$. Note that there is no easy exact
relationship between the eigenvalues of $G_n$ for various $n$ since the
eigenvectors play a role in determining any such identity. In fact, the
eigenvalues of $G_n$ and $G_{n+1}$ need not be interlaced. However, one
can consider linear eigenvalue statistics for the graph $G(n, 2d)$.
That is, for any $d$-regular graph on $n$ vertices $G$
and function $f\dvtx \mathbb{R}\rightarrow\mathbb{R}$, define the
random variable
\[
\tr f(G):=\sum_{i=1}^n \hat{f}(
\lambda_i),
\]
where $\lambda_1\ge\cdots\ge\lambda_n$\vspace*{1pt} are the eigenvalues of
adjacency matrix of $G$
\textit{divided} by $2(2d-1)^{1/2}$, and
$\hat{f}$ is $f$ with its constant term adjusted [see (\ref{naujas})
for the full definition].
The scaling is necessary to take a
limit with respect to $d$.

By a polynomial basis we refer to a sequence of polynomials $\{
f_0\equiv1, f_1,\break  f_2,  \ldots\}$ such that $f_k$ is a polynomial of
degree $k$ of a single argument over reals. In the statement below
$[\infty]$ will refer to $\NN$.\vadjust{\goodbreak}
%
%
\begin{thmm}\label{mainthm:eigenvalues}
There exists a polynomial basis $\{ f_i, i\in\NN\}$
(depending on~$d$)
such that for any $K\in\NN\cup\{ \infty\}$, the process
$( \tr f_k(G(s+t)), k \in[K], t \ge0 )$ converges in law, as
$s$ tends to infinity, to the Markov process $(N_k(t), k\in[K], t
\ge0 )$ of Theorem \ref{mainthm:cycles}.
[The polynomials are given explicitly in \eqref{eq:fbasis}.]
Hence, for any polynomial $f$, the process $ (\tr f(G(s+t)) )$
converges to a linear combination of the coordinate processes
of $(N_k(t), k\in\NN)$.
\end{thmm}

The Markov property is especially intriguing since, to the best of our
knowledge, no similar property of eigenvalues of the standard Random
Matrix ensembles is known. For the special case of minors of the
Gaussian Unitary/Orthogonal Ensembles, the entire distribution of
eigenvalues across minors of various sizes do satisfy a Markov
property. However, this is facilitated by the known symmetry properties
of the eigenvectors, and do not extend to other examples of Wigner matrices.

For our final result, we will take $d$ to infinity. We will make the
following notational convention: for any polynomial $f$, we will denote
the limiting
process of $( \tr f(G(s+t)), t \ge0 )$
by $(\tr f (G(\infty+ t) ), t\ge0 )$.
Recall that this process is a linear combination of $(N_k(t), k \in
\NN, t\ge0)$. 

%
%
\begin{thmm}\label{mainthm:chebycov}
Let $\{T_k, k \in\NN\}$ denote the Chebyshev orthogonal polynomials
of the first kind on $[-1,1]$.
As $d$ tends to infinity, the collection of processes
\[
\bigl( \tr T_{k} \bigl(G(\infty+ t) \bigr) - \E\tr T_{k}
\bigl(G(\infty+ t) \bigr), t\ge0, k\in\NN \bigr)
\]
converges weakly in $D^\infty[0,\infty)$ to a collection of independent
Ornstein--Uhlen\-beck processes $ ( U_k(t), t\ge0, k\in\NN
)$, running in equilibrium. Here the equilibrium distribution of $U_k$
is $N(0, k/2)$ and $U_k$ satisfies the stochastic differential equation
\[
d U_k(t)= - k U_k(t) \,d t + k \,d W_k(t),\qquad t
\ge0,
\]
and $(W_k, k \in\NN)$ are i.i.d. standard one-dimensional Brownian motions.

Thus, the collection of random variables
$ (\tr T_{k} (G(\infty+ t) ) -\break  \E\tr T_{k}
(G(\infty+ t) ) )$,
indexed by $k$ and $t$, converges as $d$ tends to infinity
to a centered Gaussian process with covariance kernel given by
%
%
\begin{equation}
\label{eq:covcheby} \lim_{d\rightarrow\infty} \cov\bigl( \tr T_i
\bigl(G(\infty+ t) \bigr), \tr T_k \bigl(G(\infty+ s) \bigr) \bigr)=
\delta_{ik} \frac
{k}{2}e^{k(s-t)}
\end{equation}
for $s\leq t$.
\end{thmm}

A comparison of \eqref{eq:covcheby} with Borodin's result \eqref
{eq:chebycovBor} shows that the above limit captures a key property of
the GFF covariance structure. The appearance of the exponential is
merely due to a deterministic time-change of the process. A somewhat
more detailed discussion can be found in the following section.

%
%
\begin{rmk}
A common model for random regular graphs is the \emph{configuration
model} or
\emph{pairing model} (see \cite{W} for more information). The model is
defined as follows: Start with $n$ buckets, each containing $d$ prevertices.
Then, separate these $dn$ prevertices into pairs, choosing uniformly
from every possible pairing. Finally, collapse each bucket into a
single vertex,
making an edge between one vertex and another if a prevertex in one
bucket is paired
with a prevertex in the other bucket. This model has the advantage
that choosing
a graph from it conditional on it containing no loops or parallel
edges is the same
as choosing a graph uniformly from the set of graphs without loops and
parallel edges.
The model also allows for graphs of odd degrees, unlike the
permutation model.

It is possible to construct a process of growing random regular graphs similar
to the one in this paper using
a dynamic version of this model. Given some initial pairing of prevertices
labeled $\{1,\ldots,dn\}$,
extend it to a random pairing of $\{1,\ldots,dn+2\}$ by the following
procedure:
Choose $X$ uniformly from $\{1,\ldots,dn+1\}$. Pair $dn+2$ with $X$.
If $X=dn+1$, leave the other pairs unchanged; if not, pair the
previous partner of
$X$ with $dn+1$. This is an analogue of the CRP in the setting of
random pairings,
in that if the initial pairing is uniformly chosen, then so is the
extended one.

If $d$ is odd, we repeat this procedure a total of $d$ times to extend
a random
$d$-regular graph on $n$ vertices to have
$n+2$ vertices (when $d$ is odd, the number of vertices in the graph
must be even).
When $d$ is even, repeat $d/2$ times to add one new vertex to a random graph.
In this way, we can construct a sequence of growing random regular graphs.
We believe that all the results of this paper hold in this model with
minor changes,
with similar proofs.
\end{rmk}

\subsection{Existing literature}

The study of the spectral properties of sparse regular random graphs is
motivated by several different problems. These matrices do not fall
within the purview of the standard techniques of Random Matrix Theory
(RMT) due to their sparsity and lack of independence between entries.
However, extensive simulations \cite{JMRR} point to conjectures that
these matrices still belong to the universality class of random
matrices. For example, it is conjectured via simulations \cite
{miller08} that the distribution of the second largest eigenvalue (in
absolute value) is given by the Tracy--Widom distribution. In the
physics literature, eigenvalues of random regular graphs have been
considered as a toy model of quantum chaos 
\cite{Smi10,OGS09,OS10}. Simulations suggest that
the eigenvalue spacing distribution has the same limit as that of the
Wigner matrices. A limiting \textit{Gaussian wave} character of
eigenvectors have also been conjectured \cite{Elon08,Elon09,ES10}.
Some fine properties of eigenvalues and eigenvectors can indeed be
proved for a single permutation matrix; see \cite{W00} and \cite{BAD}.

Somewhat complicating the matter is the fact that when the degree $d$
is kept fixed and we let $n$ go to infinity,
several classical results about random matrix ensembles fail. A bit
more elaboration on this point is needed. The two parameters in the
ensemble of random graphs are the degree $d$ and the order $n$. In the
permutation model it is possible to construct random regular graphs for
every possible value of $(d, n)$ where $d$ is an even positive integer
and $n$ is any positive integer. Hence, one can consider various kinds
of limits of these parameters. We will refer as the \textit{diagonal}
limit the procedure of having a sequence of $(d, n)$ where both these
parameters simultaneously go to infinity. To maintain sparsity,\setcounter{footnote}{1}\footnote{The nonsparse can be typically absorbed within standard techniques of
RMT by comparing with a corresponding Erd\H{o}s--R\'enyi graph whose
adjacency matrix has independent entries.} it is usually assumed that
$d$ is at most poly-logarithmic in $n$. No lower bound on the growth
rate of $d$ is assumed. However, results are often easier to prove when
$d$ is kept fixed and we let $n$ go to infinity. Suppose for each $d$
one gets a limiting object (say a probability distribution); one can
now take $d$ to infinity and explore limits of the sequence of these
objects. We will refer to this procedure ($\lim_{d\rightarrow\infty}
\lim_{n\rightarrow\infty}$) as the triangular limit. The triangular
limit is often identical to the diagonal limit irrespective of the
sequence through which the diagonal limit is taken, while maintaining sparsity.
Moreover, these limiting statistics frequently match with those of the
GOE ensemble and the real symmetric Wigner matrices. This is true, for
example, for the empirical spectral distribution \cite{DP,TVW} and
fluctuations of smooth linear eigenvalue statistics \cite{DJPP}.

Our present result is a triangular limit result. Let us first explain
the connection with the massless GFF. We follow Definition 2.12 and the
first example in Section~2.5 of \cite{SS07}. Consider the space of
smooth real functions compactly supported on $\mathbb{H}$ with the
Dirichlet inner product $\langle f,g \rangle=\int_{\mathbb{H}}
\nabla f\cdot
\nabla g\, dz$.
Let $H$ be the completion of this pre-Hilbert space.
The GFF can be thought of as a
random distribution $h$ which associates with every $f\in H$ a mean
zero Gaussian random variable $\langle h,f \rangle$ that is an
$\mathbf{L}^2$
isometry in the sense that $\operatorname{Cov} ( \langle h,f
\rangle,
\langle h,g \rangle )= \langle f,g \rangle$. Now, one can perform the
integration
of a function $f$ with $h$ by first integrating their traces over
semicircular arcs of a fixed radius, and then a further integral over
the radius. Over the semicircular arcs Fourier transforms (or Chebyshev
Polynomials, for real functions) provide an orthogonal basis for this
Gaussian field. As one parametrizes the radius properly, one obtains
independent Ornstein--Uhlenbeck processes for each Chebyshev polynomial.
Hence, these OU processes completely determine the GFF covariance
structure. This explains the word ``equivalent'' on page 2, paragraph 4,
and is the essence of the calculations done in \cite{Bor1}. See also~\cite{Spohn98} for a similar formalism for Dyson's Brownian motion on
the circle.

One of the reasons why we cannot prove a full GFF convergence is that
the parameters $d$ and $n$ behave independently of one another. The
degree $d$ determines the support of the spectral distribution
$[-2\sqrt{2d-1}, 2 \sqrt{2d-1}]$, asymptotically independent of $n$.
For Wigner
matrices, the dimension itself determines the length of the spectral
support. This results in the parametrization of \eqref{eq:heightfn}. It
should be possible to extend our results to a GFF convergence by either
letting $d$ grow with $n$ in the graph, or, even by letting $d$ grow
with time for the limiting Poisson structure in Theorem \ref
{mainthm:cycles}.
Though we have not attempted this
in the present article, we prove a result along these
lines in~\cite{naujasRE}.

\section{Preliminaries}

\subsection{A primer on the Chinese Restaurant Process}

The CRP, introduced by Dubins and Pitman, is a particular example of a
two parameter family of stochastic processes that constructs
sequentially random exchangeable partitions of the positive integers
via the cyclic decomposition of a random permutation. Our short
description is taken from \cite{Pit}, Section~3.1.

An initially empty restaurant has an unlimited number of circular
tables numbered $1,2,\ldots,$ each capable of seating an unlimited
number of customers. Customers numbered $1,2,\ldots$ arrive one by one
and are seated at the tables according to the following plan. Person
$1$ sits at table $1$. For $n \ge1$ suppose that $n$ customers have
already entered the restaurant, and are seated in some arrangement,
with at least one customer at each of the tables $j$ for $1\le j \le k$
(say), where $k$ is the number of tables occupied by the first $n$
customers to arrive. Let customer $n + 1$ choose with equal probability
to sit at any of the following
$n + 1$ places: to the left of customer $j$ for some $1\le j \le n$, or
alone at table $k+1$.
Define $\pi^{(n)}\dvtx [n]\to[n]$ as the permutation
whose cyclic decomposition is given by the tables; that is, if
after $n$ customers have entered the restaurant, customers $i$ and $j$ are
seated at the same table, with $i$ to the left of
$j$, then $\pi^{(n)}(i)=j$, and if customer $i$ is seated alone at
some table
then $\pi^{(n)}(i)=i$. The sequence $(\pi^{(n)})$ then has features (i)
and (ii) mentioned in the first paragraph of Section~\ref{sec:mainresults}.

\subsection{Combinatorics on words}\label{sec:wordcombinatorics}
The graph $G_n$, formed from the independent permutations
$\pi^{(n)}_1,\ldots,\pi^{(n)}_d$, can be considered as a directed,
edge-labeled
graph in a natural way. For convenience, drop superscripts
and let $\pi_l=\pi^{(n)}_l$.
If $\pi_l(i)=j$, then
by definition $G_n$ contains an edge between $i$ to $j$. When convenient,
we consider this edge to be directed from $i$ to $j$ and to be
labeled by $\pi_l$.

Consider a walk on $G_n$, viewed in this way, and imagine writing
down the label of each edge as it is traversed, putting
$\pi_i$ or $\pi_i^{-1}$ according to the direction
we walk over the edge.\vadjust{\goodbreak}
We call a walk \emph{closed} if it starts and ends at the same vertex,
and we call a closed walk a cycle
it never visits a vertex twice (besides the first and last one),
and it never traverses an
edge more than once in either direction.
Thus the word $w=w_1\cdots w_k$ formed as a cycle is traversed
is cyclically reduced, that is, $w_i\neq w_{i+1}^{-1}$ for all $i$,
considering $i$ modulo
$k$. For example,
following an edge and then immediately
backtracking does not form a $2$-cycle,
and the word formed by this walk is
$\pi_i\pi_i^{-1}$ or $\pi_i^{-1}\pi_i$ for some
$i$, which is not cyclically reduced. We consider two cycles
equivalent if they are both walks on an identical set of edges;
that is, we ignore the starting vertex and the direction of the walk.

%
\begin{figure}[b]

\includegraphics{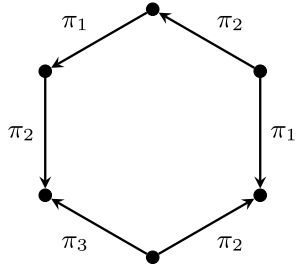}

\caption{A cycle whose word
is the equivalence class of
$\pi_2\pi_1^{-1}\pi_2\pi_1\pi_2\pi_3^{-1}$ in $\Ww_6/D_{12}$.}
\label{fig:cycleword}
\end{figure}

Let $\Ww_k$ denote the set
of cyclically reduced words of length $k$.
We would like to associate each $k$-cycle in $G_n$
with the word in $\Ww_k$ formed
by the above procedure, but since we can start the walk
at any point in the cycle and walk in either of two directions, there
are actually up to $2k$ different words that could be
formed by it. Thus, we identify
elements of $\Ww_k$ that differ only by rotation and inversion
(e.g., $\pi_1\pi_2^{-1}\pi_1\pi_2$ and
$\pi_1^{-1}\pi_2\pi_1^{-1}\pi_2^{-1}$) and denote the resulting
set by $\Ww_k/D_{2k}$, where $D_{2k}$ is the dihedral group acting
on the set $\Ww_k$ in the natural way.
%
%
\begin{defn}[(Properties of words)]
For any $k$-cycle
in $G_n$, the element of $\Ww_k/D_{2k}$ given by walking
around the cycle is called the \emph{word} of the cycle (see
Figure~\ref{fig:cycleword}).
For any word $w$, let $\vert w\vert$ denote the length of $w$.
Let $h(w)$ be the largest number $m$ such that
$w=u^m$ for some word $u$. If $h(w)=1$, we call $w$ \emph{primitive}.
For any $w\in\Ww_k$, the orbit of $w$ under the action
of $D_{2k}$ contains $2k/h(w)$ elements, a~fact which we will
frequently use.
Let $c(w)$ denote the number of pairs of double letters in $w$, that is,
the number of integers $i$ modulo $|w|$ such that $w_i=w_{i+1}$.
If~$w$ has length 1, then we define $c(w)=0$.
For example, $c(\pi_1\pi_1\pi_2^{-1}\pi_2^{-1}\pi_1)=3$.
We will also consider $|\cdot|$, $h(\cdot)$, and $c(\cdot)$ as
functions on $\Ww_k/D_{2k}$, since they are invariant
under cyclic rotation and inversion.
\end{defn}

To more easily refer to words in $\Ww_k/D_{2k}$, choose some
canonical representative $w_1\cdots w_k\in\Ww_k$ for
every $w\in\Ww_k/D_{2k}$. Based on this, we will often
think of elements of $\Ww_k/D_{2k}$ as words instead of equivalence
classes, and we will make statements about the $i$th letter
of a word in $\Ww_k/D_{2k}$. For $w=w_1\cdots w_k\in\Ww_k/D_{2k}$,
let $w^{(i)}$ refer
to the word in $\Ww_{k+1}/D_{2k+2}$ given by
$w_1\cdots w_i w_i w_{i+1}\cdots w_k$. We refer to this
operation as \emph{doubling} the $i$th letter of $w$.
A related operation is to
\emph{halve} a pair of double letters, for example producing
$\pi_1\pi_2\pi_3\pi_4$ from $\pi_1\pi_2\pi_3\pi_4\pi_1$.
(Since we apply these operations
to words identified with their rotations, we do not need
to be specific about which letter of the pair is deleted.)
The following technical lemma underpins most
of our combinatorial calculations.
%
%
\begin{lemma}\label{lem:doublehalve}
Let $u\in\Ww_k/D_{2k}$ and $w\in\Ww_{k+1}/D_{2k+2}$.
Suppose that $a$ letters in $u$ can be doubled to form $w$,
and $b$ pairs of double letters in $w$ can be halved to form $u$.
Then
\[
\frac{a}{h(u)}=\frac{b}{h(w)}.
\]
\end{lemma}
%
%
\begin{rmk}
At first glance, one might expect that $a=b$. The example
$u=\pi_1\pi_2\pi_1\pi_1\pi_2$ and $w=\pi_1\pi_1\pi_2\pi_1\pi
_1\pi_2$
shows that this is wrong, since only one letter in $u$ can
be doubled to give $w$, but two different pairs in $w$ can
be halved to give $u$.
\end{rmk}
\begin{pf}
Let $\orb(u)$ and $\orb(w)$ denote the orbits of $u$ and $w$
under the action of the dihedral group
in $\Ww_k$ and $\Ww_{k+1}$,
respectively.
When we speak of halving
a pair of letters in a word in $\orb(w)$, always delete the second of
the two letters (e.g., $\pi_1\pi_2\pi_1$ becomes $\pi_1\pi
_2$, not
$\pi_2\pi_1$). When we double a letter in a word in $\orb(u)$,
put the new letter after the doubled letter (e.g., doubling
the second letter of $\pi_1\pi_2^{-1}$ gives $\pi_1\pi_2^{-1}\pi_2^{-1}$,
not $\pi_2^{-1}\pi_1\pi_2^{-1}$).

For each of the $2k/h(u)$ words in $\orb(u)$, there are
$a$ doubling operations yielding a word in $\orb(w)$.
For each of the $(2k+2)/h(w)$ words in $\orb(w)$, there
are $b$ halving operations yielding a word in $\orb(u)$.
For every halving operation on a word in $\orb(w)$,
there is a corresponding doubling operation on a word in $\orb(u)$
and vice versa, except for halving operations that
straddle the ends of the word, as in $\pi_1\pi_2\pi_1$.
There are $2b/h(w)$ of these, giving us
\begin{eqnarray*}
\frac{2ka}{h(u)} &=& \frac{(2k+2)b}{h(w)}-\frac{2b}{h(w)}
\\
&=&\frac{2kb}{h(w)},
\end{eqnarray*}
and the lemma follows from this.
\end{pf}

Let $\Ww'=\bigcup_{k=1}^{\infty}\Ww_k/D_{2k}$, and let
$\Ww'_K=\bigcup_{k=1}^K\Ww_k/D_{2k}$.
We will use the previous lemma to prove the following
technical property of the $c(\cdot)$ statistic.
%
%
\begin{lemma}\label{lem:wordcounts}
In the vector space with basis $\{q_w\}_{w\in\Ww_K'}$,
\[
\sum_{w\in\Ww_{K-1}'}\sum_{i=1}^{|w|}
\frac{1}{h(w)}q_{w^{(i)}} = \sum_{w\in\Ww_K'}
\frac{c(w)}{h(w)}q_w.
\]
\end{lemma}
\begin{pf}
Fix some $w\in\Ww_k/D_{2k}$, and
let $a(u)$ denote the number of letters of $u$ that can be doubled to give
$w$, for any $u\in\Ww_{k-1}/D_{2k-2}$. We need to prove that
\[
\sum_{u\in\Ww_{k-1}/D_{2k-2}}\frac{a(u)}{h(u)}=\frac{c(w)}{h(w)}.
\]
Let $b(u)$ be the number of pairs in $w$ that can be halved to give $u$.
By Lemma~\ref{lem:doublehalve},
\[
\sum_{u\in\Ww_{k-1}/D_{2k-2}}\frac{a(u)}{h(u)} = \sum
_{u\in\Ww_{k-1}/
D_{2k-2}} \frac{b(u)}{h(w)},
\]
and $\sum_{u\in\Ww_{k-1}/D_{2k-2}}b(u)=c(w)$.
\end{pf}

\begin{figure}[b]

\includegraphics{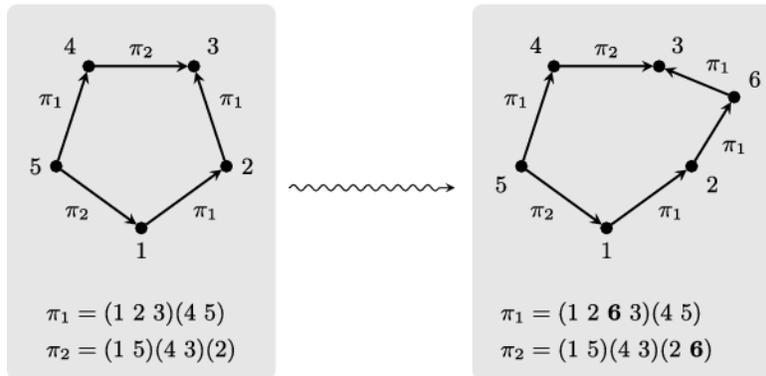}

\caption{The vertex $6$ is inserted between
vertices $2$ and $3$ in $\pi_1$, causing the above
cycle to grow.}\label{fig:growth}
\end{figure}

\section{The process limit of the cycle structure}\label{sec:lim}
\begin{figure}
\includegraphics{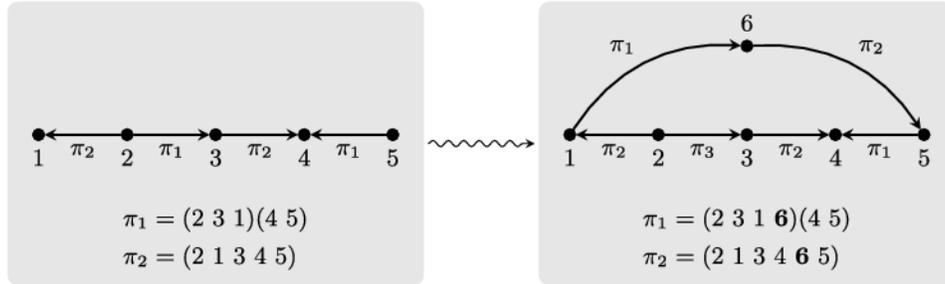}

\caption{A cycle forms ``spontaneously'' when the vertex $6$
is inserted into the graph.}\label{fig:spontaneous}
\end{figure}

As the graph $G(t)$ grows, new cycles form, which we can classify into two
types. Suppose a new vertex numbered $n$ is inserted at time
$t$, and this insertion creates a new cycle.
If the edges entering and leaving vertex
$n$ in the new cycle have the same edge label, then the new cycle has
``grown'' from a cycle with one fewer vertex, as in
Figure~\ref{fig:growth}. If the edges entering and leaving $n$ in the cycle
have different labels, then the cycle has formed ``spontaneously'' as in
Figure~\ref{fig:spontaneous}, rather
than growing from a smaller cycle.
This classification will prove essential in understanding
the evolution of cycles in $G(t)$.

\begin{figure}[b]

\includegraphics{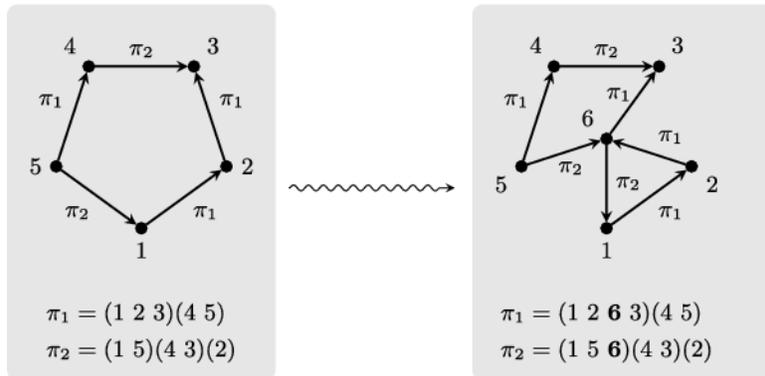}

\caption{The vertex $6$ is inserted into the cycle
in two different places in the same step, causing
the cycle to split in two. Note that
each new cycle would be classified as spontaneously
formed.}\label{fig:split}
\end{figure}

Once a cycle comes into existence in $G(t)$, it remains until a new vertex
is inserted into one of its edges. Typically, this results in the cycle
growing to a larger cycle, as in Figure~\ref{fig:growth}.
If a new vertex is simultaneously inserted into multiple
edges of the same cycle, the cycle is instead split into
smaller cycles
as in Figure~\ref{fig:split}. These new cycles
are spontaneously formed, according to the classification of new cycles
given in the previous paragraph.
Tracking the evolution of these smaller cycles in turn,
we see that as the graph evolves,
a cycle grows into a cluster of overlapping cycles.
However, it will follow from
Proposition \ref{prop:nooverlaps} that for short cycles, this behavior
is not
typical. Thus in our limiting object,
cycles will grow only into larger cycles.\label{par:split}

%

\subsection{Heuristics for the limiting process}\label{sec:heuristics}
We give some estimates that will motivate
the definition of the limiting process in Section~\ref{sec:limdef}.
This section is entirely motivational, and
we will not attempt to make anything rigorous.

Suppose that vertex $n$ is inserted into $G(t)$ at some
time $t$.
First, we consider the rate that cycles form spontaneously
with some word $w\in\Ww_k/D_{2k}$. There are
$2k/h(w)$ words in the orbit of $w$ under
the action of $D_{2k}$, and out of these, $2(k-c(w))/h(w)$
have nonequal first and last letters.
For each such word $u=u_1\cdots u_k$, we can
give a walk on the graph by starting at vertex $n$
and following the edges indicated by $u$, going from
$n$ to $u_1(n)$ to $u_2(u_1(n))$ and so on.
If this walk happens to be a cycle, the condition
$u_1\neq u_k$ implies that
it would be spontaneously formed.

In a short interval $\Delta t$ when $G(t)$ has $n-1$ vertices,
the probability that vertex $n$ is inserted is about $n \Delta t$.
For any word $u$,
the walk from vertex $n$ generated by $u$
is a cycle with probability approximately $1/n$, since
after applying the random permutations $u_1,\ldots,u_k$ in turn,
we will be left at an approximately uniform random vertex.
Any new spontaneous cycle formed with word $w$
will be counted by one of these walks,
with $u$ in the orbit of $w$, and it will be counted again
by the walk generated by
$u_k^{-1}\cdots u_1^{-1}$.
The expected number of spontaneous
cycles formed in a short interval $\Delta t$ is then approximately
\[
\frac{1}{h(w)} \bigl(k-c(w) \bigr)\frac{n \Delta t}{n}= \frac{1}{h(w)}
\bigl(k-c(w) \bigr) \Delta t.
\]
Thus, we will model the spontaneous formation of cycles with word $w$
by a Poisson process with rate $(k-c(w))/h(w)$.

Next, we consider how often a cycle with word $w\in\Ww_k$ grows into
a larger
cycle.
Suppose that $G(t)$ has $n-1$ vertices, and that it contains a cycle of
the form\vspace*{6pt}

\includegraphics{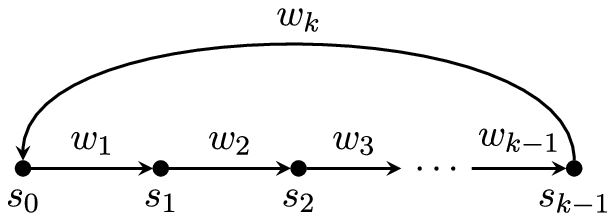}
\vspace*{6pt}

\noindent When vertex $n$ is inserted into the graph,
the probability that it is inserted after $s_{i-1}$ in permutation
$w_i$ is $1/n$. Thus, after a spontaneous cycle with word $w$
has formed, we can model the evolution of its word as a
continuous-time Markov chain
where each letter is doubled with rate one.

\subsection{Formal definition of the limiting process}\label{sec:limdef}
%
%

Consider the measure $\mu$
on $\Ww'$ given by
\[
\mu(w)= \frac{|w|-c(w)}{h(w)}.
\]
Consider a Poisson point process $\chi$ on
$\Ww' \times[0,\infty)$ with an intensity measure given by the
product measure
$\mu\otimes\leb$, where $\leb$ refers to the Lebesgue measure. Each atom
$(w,t)$ of $\chi$ represents a new spontaneous cycle with word $w$
formed at time $t$.

Now, we define
a continuous-time Markov chain on the countable space
$\Ww'$ governed by the following
rates: From state $w\in\Ww_k/D_{2k}$, jump with rate one to each
of the $k$ words
in $\Ww_{k+1}/D_{2k+2}$ obtained by doubling a letter of $w$.
If a word can be formed in more than one way by doubling a letter in $w$,
then it receives a correspondingly higher rate. For example, from
$w=\pi_1\pi_1\pi_2$, the chain jumps to $\pi_1\pi_1\pi_1\pi_2$
with rate
two and to $\pi_1\pi_1\pi_2\pi_2$ with rate one. Let $\widetilde{P}_{w}$
denote the law of this process started from $w\in\Ww'$.

Suppose we are given a realization of $\chi$. For any
atom $(w,s)$ of the countably many atoms of $\chi$,
we start an independent process $(X_{w,s}(t), t\ge0)$ with law
$\widetilde{P}_{w}$. Define the stochastic process
\[
N_{w}{({t})}:= \mathop{\sum_{(u,s)\in\chi}}_{s\leq t}
1 \bigl\{X_{u,s}(t-s)=w \bigr\}.
\]
Interpreting these processes as in the previous section,
$N_{w}{({t})}$ counts the number of cycles
formed spontaneously at time $s$ that have grown to have word $w$ at
time~$t$.

The fact that the process exists is obvious since one can define the countably
many independent Markov chains on a suitable product space.
The following lemma establishes some of its key properties.

%
\begin{lemma}\label{lem:chainprops}
Recall that $\Ww_L'=\bigcup_{k=1}^L\Ww_k/D_{2k}$.
We have the following conclusions:
\begin{longlist}[(iii)]
\item[(i)] For any $L \in\NN$, the stochastic process $\{ (N_{w}{({t})},
w\in\Ww_L'), t\ge0 \}$ is a time-homogeneous Markov
process with respect to its natural filtration, with RCLL paths.

\item[(ii)] Recall that for
$w\in\Ww_k/D_{2k}$, the element $w^{(i)}\in\Ww_{k+1}/D_{2k+2}$ is
the word formed by doubling the $i$th letter of $w$.
The generator for the Markov process
$\{ (N_{w}{({t})},
w\in\Ww_L'), t\ge0 \}$ acts on $f$ at
$x=(x_{w}, w\in\Ww_L')$ by
\begin{eqnarray*}
\mathcal{L} f(x) &= &\sum_{w\in\Ww_L'} \sum
_{i=1}^{|w|}x_{w} \bigl[
f(x-e_{w}+e_{w^{(i)}}) - f(x) \bigr]
\\
&&{}+ \sum_{w\in\Ww'_L} \frac{|w|-c(w)}{h(w)} \bigl[
f(x+e_{w}) - f(x) \bigr],
\end{eqnarray*}
where $e_{w}$ is the canonical basis vector equal to one at entry $w$ and
equal to zero everywhere else. For a word $u$ of length greater than
$L$, take $e_u=0$.

\item[(iii)] The product measure of $\Poi(1/h(w))$ over all $w\in\Ww_L'$ is
the unique invariant measure for this Markov process.
\end{longlist}
\end{lemma}

\begin{pf} Conclusion (i) follows from construction, as does conclusion
(ii). To prove conclusion (iii), we start by the fundamental identity
of the Poisson distribution: if $X\eqd\Poi(\lambda)$, then for any
function $f$, we have
%
%
%
\begin{equation}
\label{eq:poiiden} \E X g(X)= \lambda\E g(X+1).
\end{equation}

We need to show that if the coordinates of $X=(X_{w}, w\in\Ww'_L)$
are independent Poisson random variables with $\E X_{w}=1/h(w)$, then
%
%
\begin{equation}
\label{eq:poiinvpf} \E\mathcal{L}f(X)=0.
\end{equation}
Since the process is an irreducible Markov chain on countable state
space, the existence of one invariant distribution shows that the chain
is positive recurrent and that the invariant distribution is unique.

To argue \eqref{eq:poiinvpf}, we will repeatedly apply identity \eqref
{eq:poiiden} to functions $g$
constructed from $f$ by keeping all but one coordinate fixed.
Thus, for any $w\in\Ww_L'$ and $1\leq i\leq|w|$,
we condition on all $X_u$ with $u\neq w$ and hold
those coordinates of $f$ fixed to obtain,
\[
\E X_w f ( X - e_{w} + e_{w^{(i)}} )=
\frac{1}{h(w)} \E f ( X + e_{w^{(i)}} )
\]
taking $e_{w^{(i)}}=0$ when $|w|=L$.
In the same way,
\[
\E X_{w} f ( X )= \frac{1}{h(w)} \E f ( X + e_{w} ).
\]
By these two equalities,
\begin{eqnarray*}
&&\E\sum_{w\in\Ww_L'}\sum_{i=1}^{|w|} X_w
\bigl[f(X-e_w+e_{w^{(i)}})-f(X) \bigr]
\\
&&\qquad= \sum_{w\in\Ww_L'}\sum_{i=1}^{|w|}
\frac{1}{h(w)}\E \bigl[ f(X+e_{w^{(i)}}) -f(X+e_w) \bigr]
\\
&&\qquad=\sum_{w\in\Ww_{L-1}'}\sum_{i=1}^{|w|}
\frac{1}{h(w)}\E f(X+e_{w^{(i)}}) + \sum_{w\in\Ww_L/D_{2L}}
\frac{|w|}{h(w)}\E f(X)
\\
&&\qquad\quad{} - \sum_{w\in\Ww'_L}\frac{|w|}{h(w)}\E
f(X+e_w).
\end{eqnarray*}
Specializing Lemma \ref{lem:wordcounts} to $q_w=\E f(X+e_w)$,
the first sum is
\[
\sum_{w\in\Ww_{L-1}'}\sum_{i=1}^{|w|}
\frac{1}{h(w)}\E f(X+e_{w^{(i)}}) =\sum_{w\in\Ww'_L}
\frac{c(w)}{h(w)}\E f(X+e_w),
\]
which gives us
\begin{eqnarray*}
&&\E\sum_{w\in\Ww_L'}\sum_{i=1}^{|w|} X_w
\bigl[f(X-e_w+e_{w^{(i)}})-f(X) \bigr]
\\
&&\qquad= \sum_{w\in\Ww'_L}\frac{c(w)-|w|}{h(w)}\E
f(X+e_w) +\sum_{w\in\Ww_L/D_{2L}}\frac{|w|}{h(w)}
\E f(X).
\end{eqnarray*}
All that remains in proving \eqref{eq:poiinvpf} is to show that
\[
\sum_{w\in\Ww_L'}\frac{|w|-c(w)}{h(w)}=\sum
_{w\in\Ww_L/D_{2L}} \frac{|w|}{h(w)}.
\]
Specializing Lemma \ref{lem:wordcounts} to $q_w=1$ shows that
$\sum_{w\in\Ww_L'}c(w)/h(w)=\break  \sum_{w\in\Ww_{L-1}'}|w|/{h(w)}$.
Thus,
\begin{eqnarray*}
\sum_{w\in\Ww_L'}\frac{|w|-c(w)}{h(w)}&=& \sum
_{w\in\Ww_L'}\frac{|w|}{h(w)}- \sum_{w\in\Ww_{L-1}'}
\frac{|w|}{h(w)}
\\
&= &\sum_{w\in\Ww_L/D_{2L}}\frac{|w|}{h(w)},
\end{eqnarray*}
establishing \eqref{eq:poiinvpf} and completing the proof.
\end{pf}
From now on, we will consider the process $(N_{w}{({t})}, k\in\NN,
t\geq0)$
to be running under stationarity, that is, with marginal distributions
given by conclusion (iii) of the last lemma.
This process is easily constructed as described above, but with
additional point masses
of weight $1/h(w)$ for each $w\in\Ww'$ at $(w,0)$ added to the intensity
measure of $\chi$, thus giving us the correct distribution at time zero.

\subsection{Time-reversed processes}\label{sec:reversals}
Fix some time $T>0$. We define the time-reversal
$\overleftarrow{N}_{w}{({t})}:=N_{w}{({T-t})}$ for
$0\leq t\leq T$.
%
%
\begin{lemma}
For any fixed $L\in\NN$, the process
$\{(\overleftarrow{N}_{w}{({t})}, w\in\Ww'_L),\break  0\leq t\leq T\}$
is a time-homogenous Markov process with respect to the natural
filtration.
A~trivial modification at jump times renders RCLL paths.
The transition rates of this chain are given
as follows.
Let $u\in\Ww_{k-1}/D_{2k-1}$ and $w\in\Ww_k/D_{2k}$, and
suppose that $u$ can be obtained from $w$
by halving $b$ different pairs. Let $x = (x_w, w\in\Ww'_L)$.
\begin{longlist}[(iii)]
\item[(i)] The chain jumps from $x$ to
$x+e_{u}-e_{w}$ with rate $bx_{w}$.
\item[(ii)] The chain jumps from $x$ to
$x-e_{w}$ with rate $(k-c(w))x_{w}$.
\item[(iii)] If $w\in\Ww_L/D_{2L}$, then the
chain jumps from $x$ to $x+e_{w}$
with rate $L/h(w)$.
\end{longlist}
\end{lemma}
\begin{pf}
Any Markov process run backwards under stationarity is Markov.
If the chain has transition rate $r(x,y)$ from states $x$ to $y$,
then the transition rate of the backwards chain from $x$ to $y$
is $r(y,x)\nu(y)/\nu(x)$, where
$\nu$ is the stationary distribution.
We will let $\nu$ be the stationary distribution
from Lemma \ref{lem:chainprops}(iii) and
calculate the transition rates of the backwards chain,
using the rates given in
Lemma \ref{lem:chainprops}(ii).

Let $a$ denote the number of letters in $u$ that give $w$ when
doubled. The transition rate of the original chain
from $x+e_u-e_w$ to $x$ is $a(x_u+1)$, so
the transition rate of the backwards chain
from $x$ to $x+e_{u}-e_{w}$ is
\[
a(x_u+1)\frac{\nu(x+e_{k-1,c-1}-e_{k,c})}{\nu(x)} =\frac{ah(w)x_w}{h(u)},
\]
and this is equal to $bx_w$ by Lemma \ref{lem:doublehalve}.
A similar calculation shows that the transition rate from $x$ to
$x-e_{w}$ is
\[
\frac{(k-c(w))\nu(x-e_{w})}{h(w)\nu(x)} = \bigl(k-c(w) \bigr)x_w,
\]
proving (ii). The transition rate from $x$ to $x+e_{w}$
for $w\in\Ww_L/D_{2L}$
is
\[
\frac{\nu(x+e_w)}{\nu(x)}(x_{w}+1)L = \frac{L}{h(w)},
\]
which completes the proof.
\end{pf}
By definition,
\[
\overleftarrow{N}_{w}{({t})} = \mathop{\sum
_{(u,s)\in\chi}}_{s\leq
T-t} 1 \bigl\{X_{u,s}(T-t-s)=w
\bigr\}.
\]
We will modify this slightly to define the process
\[
\overleftarrow{M}_{w}{({t})}:= \mathop{ \sum
_{(u,s)\in\chi
}}_{s\leq T-t} 1 \bigl\{X_{u,s}(T-t-s)=w
\mbox{ and } \bigl|X_{u,s}(T-s)\bigr|\leq L \bigr\}.
\]
The idea is that $\overleftarrow{M}_{w}{({t})}$ is the same as
$\overleftarrow{N}_{w}{({t})}$, except
that it does not count cycles at time $t$ that had more than $L$ vertices
at time zero.
The process $(\overleftarrow{M}_{w}{({t})}, w\in\Ww_L')$
is a Markov chain with the same transition rates as
$(\overleftarrow{N}_{w}{({t})}, w\in\Ww_L')$,
except that it\vadjust{\goodbreak} does not jump from $x$ to $x+e_{w}$ for $w\in\Ww_L/D_{2L}$.
These two chains also have the same initial distribution, but
$(\overleftarrow{M}_{w}{({t})}, w\in\Ww_L')$ is not stationary (in
fact, it is eventually
absorbed at zero).

\section{Process convergence}\label{sec:processconvergence}

Recall that $C_{k}^{(s)}(t)$ is the number of cycles of length~$k$ in
the graph
$G(s+t)$, defined on page \pageref{page:chaindef}.
For $w\in\Ww'$, let $C_{w}^{(s)}(t)$ be the number of cycles in $G(s+t)$
with word $w$. We will prove that
$ (C_{w}^{(s)}(\cdot), w\in\Ww' )$ converges to a
distributional limit,
from which the convergence of $ (C_{k}^{(s)}(\cdot), k\in\NN)$
will follow.
The proof depends on knowing the limiting \emph{marginal} distribution of
$C_{w}^{(s)}(t)$. We provide this and more in the following theorem,
which should be of independent interest.

\newcommand{\I}{\mathbf{I}}
\newcommand{\Z}{\mathbf{Z}}
%
%
\begin{thmm}\label{thmm:processapprox}
Let $G_n=G(n,2d)$, a $2d$-regular random graph
on $n$ vertices from the permutation model. For any $k$,
let $\Ii_k$ be the set of all cycles of length $k$
on the complete graph $K_n$ with edge labels that form
a cyclically reduced word; these are the possible
$k$-cycles that might appear in $G_n$. Let
$\Ii=\bigcup_{k=1}^r\Ii_k$ for some integer $r$.

For any cycle $\alpha\in\Ii$, let
$I_\alpha=1\{G_n \mbox{ contains } \alpha\}$, and let $\I
=(I_\alpha,
\alpha
\in\Ii)$.
Let $\Z=(Z_\alpha, \alpha\in\Ii)$ be a vector whose coordinates
are independent Poisson random variables with $\E Z_\alpha=1/[n]_k$
for $\alpha\in\Ii_k$. Then for all $d\geq2$ and $n,r\geq1$,
\[
d_{\mathrm{TV}}(\I, \Z) \leq\frac{c(2d-1)^{2r-1}}{n} 
\]
for some absolute constant $c$, where $d_{\mathrm{TV}}(X,Y)$ denotes
the total
variation distance between the laws of $X$ and $Y$.
\end{thmm}

%
\begin{cor}\label{cor:marginals}
Let $\{Z_w, w\in\Ww'_K\}$
be a family of independent Poisson random variables
with $\E Z_w=1/h(w)$. For any fixed integer $K$ and $d\geq1$,
\begin{longlist}[(ii)]
\item[(i)] as $t\to\infty$,
\[
\bigl(C_w(t), w\in\Ww'_K \bigr) \,
{\buildrel\mathcal{L} \over
\longrightarrow}\, \bigl(Z_w, w\in\Ww'_K
\bigr);
\]
\item[(ii)] as $t\to\infty$, the probability that there exist
two
cycles of length $K$ or less sharing a vertex in $G(t)$
approaches zero.
\end{longlist}
\end{cor}
We give the proofs in the \hyperref[app]{Appendix}, along with some further
discussion.
Now, we turn to the convergence of the processes.

%
\begin{thmm}\label{thmm:processconvergence}
The process $ (C_{w}^{(s)}(\cdot), w\in\Ww' )$ converges in law
as $s\to\infty$ to
$(N_{w}{({\cdot})}, w\in\Ww')$.\label{thmm:convergence}
\end{thmm}
%
%
%
\begin{pf}
The main difficulty in
turning the intuitive ideas of Section~\ref{sec:heuristics}
into an actual proof is that $ (C_{w}^{(s)}(t), w\in\Ww' )$
is not Markov. 
We now sketch how we evade this problem.\vadjust{\goodbreak}
We will run our chain backwards, defining $\overleftarrow{G}_{s}(t)=G(s+T-t)$
for some fixed $T>0$. Then,
we ignore all of $\overleftarrow{G}_{s}(0)$ except for the subgraph
consisting of cycles
of size $L$ and smaller, which we will call ${\overleftarrow{\Gamma}{}_{s}(0)}$.
The graph ${\overleftarrow{\Gamma}_{s}(t)}$ is the evolution of this subgraph as time
runs backward, ignoring the rest of $\overleftarrow{G}_{s}(t)$.
Then, we consider the number of cycles with word $w$ in
${\overleftarrow{\Gamma}_{s}(t)}$,
which we call $\phi_w({\overleftarrow{\Gamma}_{s}(t)})$.
Choose $K \ll L$. Then $\phi_w({\overleftarrow{\Gamma}_{s}(t)})$ is likely to be the same
as $C_{w}^{(s)}(T-t)$
for any word $w$ with $|w|\leq K$.
The remarkable fact that makes $\phi_w({\overleftarrow{\Gamma}_{s}(t)})$
possible to analyze is that
if ${\overleftarrow{\Gamma}_{s}(0)}$
consists of disjoint cycles, then $ (\phi_w({\overleftarrow{\Gamma}_{s}(t)}),
w\in\Ww'_L )$
is a Markov chain governed by the same transition rates
as $ (\overleftarrow{M}_{w}{({t})}, w\in\Ww'_L )$.

Another important idea of the proof is to ignore the vertex labels
in $\overleftarrow{G}_{s}(t)$, so that we do not know in what order
the vertices
will be removed. Thus, we can view $\overleftarrow{G}_{s}(t)$ as a Markov
chain with the following description:
Assign each
vertex an independent $\Exp(1)$ clock. When the clock of vertex $v$
goes off, remove it from the graph, and patch together the $\pi_i$-labeled
edges entering and leaving $v$ for each $1\leq i\leq d$.

\begin{step}[{[Definitions of ${\overleftarrow{\Gamma}_{s}(t)}$ and $\phi_w$ and analysis of $
(\phi_{w}({\overleftarrow{\Gamma}_{s}(t)}),\break  w\in\Ww_L' )$]}]

Fix $T>0$ and define $\overleftarrow{G}_{s}(t)=G(s+T-t)$. As mentioned above,
we will consider
$\overleftarrow{G}_{s}(t)$ only up to relabeling of vertices, which
makes it
a process on the countable state space consisting
of all edge-labeled graphs on finitely many unlabeled vertices.
With respect to its natural filtration, it is a Markov chain
in which each vertex is removed with rate one,
as described above.

To formally define
${\overleftarrow{\Gamma}_{s}(t)}$, fix integers $L>K$ and
let ${\overleftarrow{\Gamma}_{s}(0)}$ be the subgraph of $\overleftarrow
{G}_{s}(0)$ made up of all
cycles of length $L$ or less.
We then evolve ${\overleftarrow{\Gamma}_{s}(t)}$ in parallel with $\overleftarrow{G}_{s}(t)$.
When a vertex $v$ is deleted from $\overleftarrow{G}_{s}(t)$, the corresponding
vertex $v$ in ${\overleftarrow{\Gamma}_{s}(t)}$ is deleted if it is present. If $v$
has a $\pi_i$-labeled edge entering and leaving it in
${\overleftarrow{\Gamma}_{s}(t)}$, then these two edges
are patched together. Other edges
in ${\overleftarrow{\Gamma}_{s}(t)}$ adjacent to $v$ are deleted.
This makes ${\overleftarrow{\Gamma}_{s}(t)}$ a subgraph of $\overleftarrow{G}_{s}(t)$,
as well as a
continuous-time Markov chain on the countable
state space consisting
of all edge-labeled graphs on finitely many unlabeled vertices.
The transition probabilities of ${\overleftarrow{\Gamma}_{s}(t)}$ do not depend on~$s$.

From Corollary \ref{cor:marginals}, we can find the limiting
distribution of ${\overleftarrow{\Gamma}_{s}(0)}$.
Suppose that $\gamma$ is a graph in the process's state space
that is not a disjoint union
of cycles. By Corollary \ref{cor:marginals}(ii),
\[
\lim_{s\to\infty}\P\bigl[ {
\overleftarrow{\Gamma}_{s}(0)}=\gamma\bigr]=0.
\]
Suppose instead that
$\gamma$ is made up of disjoint cycles,
with $z_w$ cycles of word $w$ for each $w\in\Ww'_L$.
By Corollary \ref{cor:marginals}(i),
%
%
%
\begin{equation}
\lim_{s\to\infty}\P\bigl[ {
\overleftarrow{\Gamma}_{s}(0)}=\gamma\bigr]=\prod
_{w\in\Ww_L'} \P[Z_w=z_w],
\label{eq:Gammadist}\vadjust{\goodbreak}
\end{equation}
where $(Z_w, w\in\Ww'_L)$ are independent Poisson random
variables with $\E Z_w=1/h(w)$.
Thus, ${\overleftarrow{\Gamma}_{s}(0)}$ converges in law as $s\to\infty$ to
a limiting distribution supported on the graphs
made up of disjoint unions of cycles.
For different values of~$s$, the chains ${\overleftarrow{\Gamma}_{s}(t)}$
differ only in their
initial distributions, and the convergence in law
of ${\overleftarrow{\Gamma}_{s}(0)}$ as $s\to\infty$ induces the
process convergence
of $\{{\overleftarrow{\Gamma}_{s}(t)}, 0\leq t\leq T\}$
to a Markov
chain $\{{\overleftarrow{\Gamma}_{}(t)}, 0\leq t\leq T\}$ with the same
transition rates
whose initial distribution is the limit
of ${\overleftarrow{\Gamma}_{s}(0)}$.

For any finite
edge-labeled graph $G$, let $\phi_{w}(G)$ be the number of
cycles in $G$ with word $w$.
By the continuous mapping theorem, the process
$(\phi_{w}({\overleftarrow{\Gamma}_{s}(t)}),
w\in\Ww_L')$ converges in law to $(\phi_{w}({\overleftarrow{\Gamma}_{}(t)}), w\in
\Ww_L')$
as $s\to\infty$.

We will now demonstrate that this process has the same law
as $(\overleftarrow{M}_{w}{({t})}, w\in\Ww_L')$.
The graph ${\overleftarrow{\Gamma}_{}(t)}$ consists of disjoint cycles at time $t=0$,
and as it evolves, these cycles shrink or are destroyed. The process
$(\phi_{w}({\overleftarrow{\Gamma}_{}(t)}), w\in\Ww_L')$
jumps exactly when a vertex in a cycle in ${\overleftarrow{\Gamma}_{}(t)}$ is deleted.
If the deleted
vertex lies in a cycle
between two edges with the same label, the cycle shrinks.
If the deleted vertex lies in a cycle between two edges with different
labels, the cycle is destroyed.
The only relevant consideration in where the process
will jump at time $t$ is the number of vertices of these
two types in ${\overleftarrow{\Gamma}_{}(t)}$, which can be deduced
from
$(\phi_{w}({\overleftarrow{\Gamma}_{}(t)}), w\in\Ww_L')$.
Thus, this process is a Markov chain.

Consider two words $u,w\in\Ww_K'$ such that $w$ can be
obtained from $u$ by
doubling a letter. Suppose that $u$ can be obtained from
$w$ by halving any of $b$ pairs of letters. Suppose that
the chain is at state
$x=(x_v, v\in\Ww_L')$. There are $bx_w$ vertices that when
deleted cause the chain to jump from $x$ to $x-e_w+e_u$, each
of which is removed with rate one. Thus, the chain jumps
from $x$ to $x-e_w+e_u$ with rate $bx_w$. Similarly, it
jumps to $x-e_w$ with rate $(|w|-c(w))x_w$. These are
the same rates as the chain $(\overleftarrow{M}_{w}{({t})}, w\in\Ww_L')$
from Section~\ref{sec:reversals}. The initial distribution
given by \eqref{eq:Gammadist} is also the same as that
of $(\overleftarrow{M}_{w}{({t})}, w\in\Ww_L')$, demonstrating that
the two processes
$ (\phi_{w}({\overleftarrow{\Gamma}_{}(t)}), w\in\Ww_L' )$
and $(\overleftarrow{M}_{w}{({t})}, w\in\Ww_L')$ have the same law.
\end{step}

\begin{step}[{[Approximation of $\overleftarrow{C}{}^{(s)}_{w}(t)$ by $\phi
_{w}({\overleftarrow{\Gamma}_{s}(t)})$]}]

We will compare the two processes
$\{(\overleftarrow{C}{}^{(s)}_{w}(t), w\in\Ww'_K), 0\leq t\leq T \}$
and $\{(\phi_w({\overleftarrow{\Gamma}_{s}(t)}), w\in\Ww'_K), 0\leq t\leq T\}$
and
show that for sufficiently large $L$, they are identical with
probability arbitrarily close to one.

Consider some cycle in $\overleftarrow{G}_{s}(t)$; we can divide its vertices
into those that lie between two edges of the cycle with different labels,
and those that lie between two edges with the same label.
We call this second class the \emph{shrinking vertices} of the cycle,
because if one is deleted from $\overleftarrow{G}_{s}(t)$ as it
evolves, the cycle
shrinks. We define $E_s(L)$ to be the event that for some cycle
in $\overleftarrow{G}_{s}(0)$ of size $l>L$, at least $l-K$ of its
shrinking vertices
are deleted by time $T$.

We claim that outside of the event $E_s(L)$, the two processes
$\{(\overleftarrow{C}{}^{(s)}_{w}(t), w\in\Ww'_K), 0\leq t\leq T\}$
and $\{(\phi_w({\overleftarrow{\Gamma}_{s}(t)}), w\in\Ww'_K), 0\leq t\leq T\}$
are identical. Suppose that these two processes are not identical.
Then there is some cycle $\alpha$ of size $K$ or less present in
$\overleftarrow{G}_{s}(t)$ but not in ${\overleftarrow{\Gamma}_{s}(t)}$ for $0<t\leq T$.
As explained in Section~\ref{sec:lim}, as a cycle evolves (in forward
time), it grows into an overlapping cluster of cycles.
Thus, $\overleftarrow{G}_{s}(0)$ contains some cluster of overlapping cycles
that shrinks to $\alpha$ at time~$t$. One of the cycles in this cluster
has length greater than $L$, or the cluster would be contained
in ${\overleftarrow{\Gamma}_{s}(0)}$ and $\alpha$ would have been
contained in ${\overleftarrow{\Gamma}_{s}(t)}$.

To see that $l-K$ shrinking vertices must be deleted from this
cycle, consider the evolution of $\alpha$ into the cluster
of cycles in both forward and reverse time. If a vertex is inserted into
a single edge of a cycle in forward time, we see in reverse time the deletion
of a shrinking vertex. If a vertex is simultaneously inserted into
two edges of a cycle, causing the cycle to split,
we see in reverse time the deletion of a nonshrinking
vertex of a cycle. As $\alpha$ grows, a cycle of size greater than $L$
can form
only by single-insertion of at least $l-K$ vertices into the eventual
cycle. In reverse time, this is seen as deletion of $l-K$ shrinking
vertices.
This demonstrates
that $E_s(L)$ holds.

We will now show that for any $\eps>0$, there is an $L$ sufficiently large
that
$\P[E_s(L)]<\eps$ for any $s$.
Let $w\in\Ww_l/D_{2l}$ with $l>L$, and let
$I\subseteq[l]$
such that $\vert I\vert=l-K$ and
$w_i=w_{i-1}$ for all $i\in I$,
considering indices modulo $l$.
For any cycle in $\overleftarrow{G}_{s}(0)$
with word $l$, the set $I$ corresponds to a set of $l-K$
shrinking vertices of the cycle.

We define
$F(w,I)$ to be the event that $\overleftarrow{G}_{s}(0)$ contains one
or more
cycles with word $w$,
and that the vertices corresponding to $I$ in one of these
cycles are all deleted within
time $T$.
By a union bound,
%
%
%
\begin{equation}
\P \bigl[E_s(L) \bigr]\leq\sum_{w,I}\P
\bigl[F(w,I) \bigr].\label{eq:EsL}
\end{equation}

We proceed by enumerating all pairs of $w$ and $I$.
For any pair $w, I$, deleting the letters in $w$ at positions given
by $I$ results in a word $u\in\Ww_K/D_{2K}$.
For any given $u=u_1\cdots u_K\in\Ww_K/D_{2K}$, the word $w\in\Ww_l/D_{2l}$
must have the form
\[
w= \underbrace{u_1\cdots u_1}_{a_1\ \mathrm{times}}
\underbrace{u_2\cdots u_2}_{a_2\ \mathrm{times}}\cdots\cdots
\underbrace{u_K\cdots u_K}_{a_K\ \mathrm{times}},
\]
with $a_i\geq1$ and $a_1+\cdots+a_K=l$.
The number of choices for $a_1,\ldots,a_K$ is ${l-1\choose K-1}$,
the number of compositions of $l$ into $K$ parts, and each
of these corresponds to a choice of $w$ and $I$. There are fewer
than $a(d,K)$
choices for $u$, giving us a bound of $a(d,K){{l-1 \choose K-1}}$
choices of pairs $w$ and $I$ for any fixed $l>L$.

Next, we will show that
for any pair $w$ and $I$ with $|w|=l$,
%
%
%
\begin{equation}
\P \bigl[F(w,I) \bigr]\leq \bigl(1-e^{-T} \bigr)^{l-K}.
\label{eq:FuI}
\end{equation}
Condition on $\overleftarrow{G}_{s}(0)$ having $n$ vertices.
Consider any of the $[n]_l$ possible sequences of $l$ vertices.
Choose some representative $w'\in\Ww_l$ of $w$.
For each of these sequences,
the probability that
it forms a cycle with word $w'$ is at most $1/[n]_l$ (recall the original
definition of our random graphs in terms of random permutations).
Given that
the sequence forms a cycle, the probability
that the vertices of the cycle at positions
$I$ are all deleted within time $T$
is $(1-e^{-T})^{l-K}$. Hence
\begin{eqnarray*}
\P \bigl[F(w,I)\mid\overleftarrow{G}_{s}(0)\mbox{ has $n$
vertices}\bigr] & \leq&[n]_l\frac{1}{[n]_l} \bigl(1-e^{-T}
\bigr)^{l-K},
\\
&\leq& \bigl(1-e^{-T} \bigr)^{l-K}.
\end{eqnarray*}
This holds for any $n$,
establishing \eqref{eq:FuI}.

Applying all of this to \eqref{eq:EsL},
\[
\P \bigl[E_s(L) \bigr] \leq\sum_{l=L+1}^{\infty}a(d,K)
\pmatrix{l-1
\cr
K-1} \bigl(1-e^{-T} \bigr)^{l-k}.
\]
This sum converges, which means that for any $\eps>0$,
we have $\P[E_s(L)]<\eps$ for large enough $L$, independent of
$s$.
\end{step}

\begin{step}[{[Approximation of $\overleftarrow{N}_{w}{({t})}$ by
$\overleftarrow{M}_{w}{({t})}$]}]

Recall that we defined the processes
$\{(\overleftarrow{M}_{w}{({t})}, w\in\Ww'_K), 0\leq t\leq T\}$
and $\{(\overleftarrow{N}_{w}{({t})}, w\in\Ww'_K), 0\leq t\leq T\}$
on the same probability space. We will show that for sufficiently
large $L$, the two processes
are identical with probability arbitrarily close to one.

By their definitions, these two processes
are identical unless one of the processes $X_{u,s}(\cdot)$ started
at each atom of $\chi$ grows from a word of size $K$ or less to
a word of size $L+1$ before time $T$; we call this event $E(L)$.
Let
\[
Y= \bigl| \bigl\{(u,s)\in\chi\dvtx |u|\leq K, s\leq T \bigr\} \bigr|,
\]
the number of processes starting from a word of size $K$ or less
before time~$T$.

Suppose that $X(\cdot)$ has law
$\widetilde{P}_{w}$ for some word $w\in\Ww_k/D_{2k}$.
We can choose $L$ large enough that
$\P[|X(T)|>L ]<\eps$ for all $k\leq K$.
Then
$\P[E(L)\mid Y]<\eps Y$ by a union bound, and so
$\P[E(L)]<\eps\E Y$. Since $\E Y<\infty$, we can
make $\P[E(L)]$ arbitrarily small by choosing sufficiently large $L$.
\end{step}

\begin{step}[{[Weak convergence of $\{(\overleftarrow{C}{}^{(s)}_{w}(t), w\in
\Ww'_K), 0\leq t\leq T \}$ to $\{(\overleftarrow{N}_{w}{({t})},\break   w\in
\Ww'_K), 0\leq t\leq T\}$]}]

If two processes are identical with probability
$1-\eps$, then the total variation distance between their laws is at most
$\eps$. Thus, by steps 2 and 3, we can choose $L$ large enough that
the laws of the processes
$\{(\overleftarrow{C}{}^{(s)}_{w}(t), w\in\Ww'_K), 0\leq t\leq T\}$ and
$\{(\phi_w(
{\overleftarrow{\Gamma}_{s}(t)}, w\in\Ww'_K),
0\leq t\leq T)\}$ are arbitrarily close in total variation distance,
uniformly in $s$, and so that the laws of
$\{(\overleftarrow{M}_{w}{({t})}, w\in\Ww'_K), 0\leq t\leq T\}$ and
$\{(\overleftarrow{N}_{w}{({t})}, w\in\Ww'_K), 0\leq t\leq T\}\}$
are arbitrarily close in total variation distance.
Since total variation distance dominates the Prokhorov metric (or any other
metric for the topology of weak convergence), we can
choose $L$ such that these two pairs are each within $\eps/3$ in
the Prokhorov metric.
Since $\{(\phi_w(
{\overleftarrow{\Gamma}_{s}(t)}), w\in\Ww'_K), 0\leq t\leq T\}$
converges in law to $\{(\overleftarrow{M}_{w}{({t})}, w\in\Ww'_K),
0\leq t\leq T\}$
as $s\to\infty$, there is an $s_0$ such that for all $s\geq s_0$, the laws
of these processes are within $\eps/3$ in the Prokhorov metric.
We have thus shown that for every $\eps>0$,
the laws of $\{(\overleftarrow{C}{}^{(s)}_{w}(t), w\in\Ww_K'), 0\leq
t\leq T\}$
and
$\{(\overleftarrow{N}_{w}{({t})}, w\in\Ww'_K), 0\leq t\leq T\}$
are within $\eps$
for sufficiently large $s$, which proves that the first random vector
converges in law to the second as $s\to\infty$.
\end{step}

\begin{step}[{[Weak convergence of $\{(C_{w}^{(s)}(t), w\in\Ww'), t\geq 0\}
$ to $\{(N_{w}{({t})}, w\in\Ww'), t\geq0\}$]}]

It follows
immediately from the previous step that
the (not time-reversed) process
$\{(C_{w}^{(s)}(t), w\in\Ww'_K), 0\leq t\leq T\}$
converges in law to
$\{(N_{w}{({t})}, w\in\Ww'_K), 0\leq t\leq T\}$ for any $T>0$.
By Theorem 16.17 in \cite{Bil},
$\{(C_{w}^{(s)}(t), w\in\Ww'_K), t\geq0\}$
converges in law to
$\{(N_{w}{({t})}, w\in\Ww'_K), t\geq0\}$, which also proves that
$\{(C_{w}^{(s)}(t), w\in\Ww'), t\geq0\}$
converges in law to
$\{(N_{w}{({t})}, w\in\Ww'),\break  t\geq0\}$.~\qed 
\end{step}
\noqed\end{pf}

\begin{pf*}{Proof of Theorem \ref{mainthm:cycles}}
We now consider the case of short cycles in the graph.
We will express these as functionals
of $ (C_{w}^{(s)}(t), w\in\Ww' )$. For example, consider the count
of cycles of size $k\in\NN$. Then
$C_{k}^{(s)}(t)=\sum_{w\in\Ww_k/D_{2k}}C_{w}^{(s)}(t)$ is the number
of $k$-cycles
in $G(s+t)$, and let
\[
N_{k}{({t})}=\sum_{w\in\Ww_k/D_{2k}}N_{w}{({t})}.
\]
It follows immediately from
the continuous mapping theorem that $\{(C_{k}^{(s)}(t),\break  k\in\NN),
t\geq
0\}$
converges in law to $\{(N_{k}{({t})}, k\in\NN), t\geq0\}$ as $s\to
\infty$.

It is not hard to see that this limit is Markov and admits the following
representation: Cycles of size $k$ appear spontaneously with rate
$\sum_{w\in\Ww_k/D_{2k}}\mu(w)$. The size of each
cycle then grows as a pure birth process with generator
$L f(i) = i ( f(i+1) - f(i) )$.
The only thing we need to verify is that
%
%
\begin{equation}
\label{eq:ratecount} \sum_{w\in\Ww_k/D_{2k}}\mu(w)=\sum
_{w\in\Ww_k/D_{2k}}\frac
{k-c(w)}{h(w)}=\bigl(a(d,k)-a(d,k-1)\bigr)/2.
\end{equation}

However, this follows from Lemma \ref{lem:wordcounts} in the following
way. From that lemma, we get
\[
\sum_{w \in\Ww_k/D_{2k}} \frac{c(w)}{h(w)}= (k-1) \sum
_{w\in\Ww
_{k-1}/D_{2(k-1)}}\frac{1}{h(w)}.
\]
Thus,
\[
\sum_{w\in\Ww_k/D_{2k}}\mu(w)=\sum_{w \in\Ww_k/D_{2k}}
\frac
{k}{h(w)} - \sum_{w\in\Ww_{k-1}/D_{2(k-1)}}\frac{k-1}{h(w)}.
\]
However, the two terms on the right-hand side of the above equation are
simply half the total number of cyclically reduced words possible, of
size $k$ and $k-1$, respectively. The total number of cyclically
reduced words of size $k$ on an alphabet of size $d$ is $a(d,k)$ (see
Appendix of \cite{DJPP}). This shows \eqref{eq:ratecount} and completes
the proof.
\end{pf*}

We end with the following corollary.

%
\begin{cor}\label{cor:cyclecov}
For any $s < t$ and $j,k\in\NN$, one has:
\begin{eqnarray*}
&&\cov \bigl(N_k(t), N_j(s) \bigr) \\
&&\qquad= \cases{
\displaystyle\frac{a(d,j)}{2j}\pmatrix{k-1
\cr
k-j} p^{j} (1-p)^{k-j},&\quad
$\mbox{$p=e^{s-t},$ if $k \ge j$},$\vspace*{2pt}
\cr
0,&\quad $\mbox{otherwise}.$}
\end{eqnarray*}
\end{cor}

\begin{pf}
We will refer to the Yule processes counted by $N_k(t)$ as cycles of length
$k$ present at time $t$,
even though these ``cycles'' in the limiting process have no connection
to graphs.
If $k < j$, every cycle that is of length $j$ at time $s$ cannot grow
to a cycle of length $k$ at time $t$. Thus, $N_k(t)$ depends on
cycles that are independent of those that make up $N_j(s)$. Hence
$N_k(t)$ is independent of $N_j(s)$.

If $k \ge j$, notice that one has the following decomposition:
%
%
\begin{equation}
\label{eq:decompN} N_k(t)= \sum_{j=1}^k
\alpha(j,k) N_j(s) + Z,
\end{equation}
where $\alpha(j,k)$ is the proportion of one-dimensional pure-birth
Yule processes that were at state $j$ at time $s$ and grew to state $k$
at time $t$, and $Z$ is a random variables that counts the number\vadjust{\goodbreak} of
new births in the time interval $(s,t)$ that grew to state $k$ at time
$t$. Note that, under our invariant distribution all random variables
$\{Z, N_j(s), 1\le j \le k\}$ are independent of one another. Thus,
our conclusion follows once we show
%
%
\begin{equation}
\label{eq:expecalpha} \E\alpha(j,k)= \pmatrix{k-1
\cr
k-j} p^{j}
(1-p)^{k-j},\qquad p=e^{s-t}.
\end{equation}

The expected proportion $\E\alpha(j,k)$ is the probability that a
one-dimensional process $X_{j,k}$, with law of an Yule process starting
at $j$, is at state $k$ at time $(t-s)$.
If $\xi_j, \ldots, \xi_k$ are independent exponential
random variables with rates $j,\ldots,k$, then
\[
\E\alpha(j,k) =\P \bigl[ \{ \xi_j + \cdots+ \xi_{k-1} \le
t-s \} \cap\{ \xi_j+ \cdots+ \xi_k > t-s \} \bigr].
\]
We now use the R\'enyi representation: suppose $Y_1, Y_2, \ldots, Y_k$
are i.i.d.
$\Exp(1)$ random variables. Define the order statistics $Y_{(1)} \ge
Y_{(2)} \ge\cdots\ge Y_{(k)}$. Then, the following equality holds in
distribution
\[
( Y_{(i)} - Y_{(i+1)}, j\le i \le k )= ( \xi_i, j
\le i \le k ).
\]
Here we have defined $Y_{(k+1)}\equiv0$. Thus, in distribution,
\[
\xi_j + \cdots+ \xi_{k-1} = Y_{(j)} -
Y_{(k)}, \xi_j+ \cdots+ \xi_k =
Y_{(j)}.
\]

Thus, 
\[
\E\alpha(j,k) = P (t-s < Y_{(j)} \leq Y_{(k)} + t-s ).
\]
%
Note that, by an elementary symmetry argument, for any $u> (t-s)$, we have
\begin{eqnarray*}
&&\P \bigl[ Y_{(j)}\in(u, u+\,\d u), Y_{(j)}- Y_{(k)}
< t-s \bigr]
\\
&&\qquad=\P \bigl[ \mbox{$Y_i=u$ for some $i$, exactly $j-1$ of
$Y_1,\ldots, Y_k$ are greater than $u$,}
\\
&&\hspace*{105pt}\qquad \mbox{and the rest of $Y_1,\ldots,Y_k$ are in
$[u-t+s, u]$} \bigr]\,\d u
\\
&&\qquad= k e^{-u} \pmatrix{k-1
\cr
j-1} e^{-(j-1)u} \bigl[
e^{-u+t-s} - e^{-u} \bigr]^{k-j} \,\d u
\\
&&\qquad=k \pmatrix{k-1
\cr
j-1} e^{-ku} \bigl( e^{t-s} - 1
\bigr)^{k-j}\,\d u.
\end{eqnarray*}
Integrating out $u$ in the interval $(t-s, \infty)$, we get
\begin{eqnarray*}
\P [t-s < Y_{(j)} < Y_{(k)} + t-s ]&=&\pmatrix{k-1
\cr
j-1}
\bigl( e^{t-s} - 1 \bigr)^{k-j} \int_{t-s}^\infty
k e^{-k u} \,\d u
\\
&=& \pmatrix{k-1
\cr
j-1} \bigl( e^{t-s} - 1 \bigr)^{k-j}
e^{-k(t-s)}\\
&=& \pmatrix{k-1
\cr
j-1} e^{j(s-t)} \bigl(1 -
e^{s-t} \bigr)^{k-j}.
\end{eqnarray*}
This shows \eqref{eq:expecalpha} and completes the proof of the
corollary.\vadjust{\goodbreak}
\end{pf}

\subsection{Two-dimensional convergence}
So far, we have considered $d$ as a constant. We now view
it as a parameter of the graph and allow it to vary.
Recall that $(\Pi_d, d\in\NN)$ are independent towers of random
permutations, with $\Pi_d=(\pi_d^{(n)}, n\in\NN)$, and that
$G(n,2d)$ is defined from $\pi_1^{(n)},\ldots,\pi_d^{(n)}$.
For each $d$, we follow the construction used to define
$G(t)$ and construct $G_{2d}(t)$, a continuous-time version
of $(G(n,2d), n\in\NN)$. Let $\Ww'(d)$ be the set of equivalence
classes of cyclically reduced words as before, with the parameter $d$
made explicit. Define $C_{d,k}^{(s)}(t)$
as the number of $k$-cycles in $G_{2d}(s+t)$ and consider
the convergence of the two-dimensional field $\{(C_{d,k}^{(s)}(t),
d,k\in\NN),\break  t\geq0\}$ as $s\to\infty$.

Again, we will consider this process as a functional of another one.
Define $\Ww'(\infty)=\bigcup_{d=1}^{\infty}\Ww'(d)$,
noting that $\Ww'(1)\subseteq\Ww'(2)\subseteq\cdots$.
For any $w\in\Ww'(d)$, the number of cycles in $G_{2d'}(s+t)$
with word $w$ is the same for all $d'\geq d$.
We define $C_{w}^{(s)}(t)$ by this, so that
\[
C_{d,k}^{(s)}(t)=\mathop{\sum_{w\in\Ww'(d)}}_{|w|=k}
C_{w}^{(s)}(t).
\]
Then we will
prove convergence of $\{(C_{w}^{(s)}(t), w\in\Ww'(\infty)), t\geq
0\}$
as $s\to\infty$.

To define a limit for this process,
we extend $\mu$ to a measure on all of $\Ww'(\infty)$ and define
the Poisson point process $\chi$ on $\Ww'(\infty)\times[0,\infty)$.
The rest of the construction is identical to the one in
Section~\ref{sec:limdef}, giving
us random variables $ (N_{w}{({t})}, w\in\Ww'(\infty) )$.

%
\begin{thmm}
The process $ (C_{w}^{(s)}(\cdot), w\in\Ww'(\infty) )$ converges
in law
as\break  $s\to\infty$ to
$ (N_{w}{({\cdot})}, w\in\Ww'(\infty) )$.
\end{thmm}
\begin{pf}
It suffices to prove that
$ (C_{w}^{(s)}(\cdot), w\in\Ww'(d) )$ converges in law
as $s\to\infty$ to
$(N_{w}{({\cdot})}, w\in\Ww'(d))$ for each $d$, which we did in
Theorem \ref{thmm:convergence}.
\end{pf}
\begin{pf*}{Proof of Theorem \ref{mainthm:intertwining}}
Let
\[
N_{d,k}{({t})}= \mathop{\sum_{w\in\Ww'(d)}}_{|w|=k}
N_{w}{({t})}.
\]
By the continuous mapping theorem, $(N_{d,k}{({\cdot})}, d,k\in\NN)$
is the limit of $(C_{d,k}^{(s)}(\cdot),\break   d,k\in\NN)$ as $s\to\infty$.

Let us now describe what the limiting process is. It is obvious that
$(N_{d,k}{({\cdot})}, k\in\NN, d \in\NN)$ is jointly Markov. For
every fixed $d$, the law of the corresponding marginal is given by
Theorem \ref{mainthm:cycles}. To understand the relationship across
$d$, notice that cycles of size $k$ for $(d+1)$ consist of cycles of
size $k$ for $d$ and the extra ones that contain an edge labeled by
$\pi
_{d+1}$ of $\pi_{d+1}^{-1}$. Thus
\[
N_{d+1, k}(t) - N_{d, k}(t)= \mathop{\sum
_{w\in\Ww'(d+1)\setminus
\Ww
'(d)}}_{ \vert w\vert=k} N_w(t).
\]
This process is independent of 
$(N_{i, \cdot}$, $i\in[d])$, since the set of words involved are
disjoint. Moreover, the rates for this process are clearly the
following: cycles of size $k$ grow at rate $k$ and new cycles of size
$k$ appear at rate $[a(d+1, k) - a(d+1, k-1) - a(d,k) + a(d, k-1)]/2$.
This completes the proof of the result.
\end{pf*}

\section{Process limit for linear eigenvalue statistics}

Let us recall some of the basic facts established in \cite{DJPP},
Sections~3 and~5, that connect linear eigenvalue statistics with cycle counts. A
closed nonbacktracking walk is a walk
that begins and ends at the same vertex, and that never follows
an edge and immediately follows that same edge backwards.
If the last step of a closed nonbacktracking walk is anything other
than the reverse of the first step, we say that the walk is \emph{cyclically
nonbacktracking}. 
Cyclically nonbacktracking walks on $G_n$
are exactly the closed nonbacktracking walks whose words
are cyclically reduced.
Let $\mathrm{CNBW}_{k}^{(n)}$ denote the number of closed cyclically
nonbacktracking walks of length $k$ on $G_n$.

Cyclically nonbacktracking walks are useful because they can
be enumerated by linear functionals of a graph's eigenvalues.
Let $\{T_n(x)\}_{n \in\mathbb{N}}$ be the Chebyshev polynomials of the
first kind
on the interval $[-1,1]$.
We define a set of polynomials
%
\begin{eqnarray*}
\Gamma_0(x) & =& 1,
\\
\Gamma_{2k}(x) &=& 2T_{2k}(x)+ \frac{2d-2}{(2d-1)^k}\qquad
\forall k \geq1,
\\
\Gamma_{2k+1}(x) &=& 2T_{2k+1}(x) \qquad  \forall k \geq0.
\end{eqnarray*}

Let $A_n$ be the adjacency matrix of $G_n$, and
let $\lambda_1\ge\cdots\ge\lambda_n$ be the eigenvalues
of $(2d-1)^{-1/2}A_n/2$. 
Then
%
%
%
\begin{equation}
\sum_{i=1}^n\Gamma_k(
\lambda_i)=(2d-1)^{-k/2}\mathrm{CNBW}_{k}^{(n)}.
\label{eq:eigencnbw}
\end{equation}
Suppose that $f(x)$ is a polynomial with the expansion
$f(x) = \sum_{j=0}^k a_j \Gamma_j(x)$. We define $\tr f(G_n)$ as
\begin{equation}\label{naujas}
\tr f(G_n) := \sum_{i=1}^n f(\lambda_i) - na_0.
\end{equation}
Subtracting off the constant term $a_0$ keeps
$\tr f(G_n)$ of constant order as $n$ grows.

Now, for any cycle in $G_n$ of length $j\divides k$, we obtain $2j$
nonbacktracking walks
of length $k$ by choosing a starting point and direction and then
walking around the cycle repeatedly.
In \cite{DJPP}, Corollary 18,
it is shown that with certain conditions on the growth
of $d$ and $r$, all cyclically nonbacktracking walks of length $r$ or less
have this form with high probability.
Thus, the random vectors $ (\mathrm{CNBW}_{k}^{(n)}, 1\leq k\leq r )$
and $ (\sum_{j\divides k}2jC_{k}^{(n)}, 1\leq k\leq r )$ have
the same
limiting distribution, and the problem of finding the limiting distributions
of polynomial linear eigenvalue statistics is reduced to finding limiting
distributions of cycle counts.
We will prove Theorem \ref{mainthm:eigenvalues} by arguing that
this holds for the entire process $(G(t), t \ge0)$.

Call a cyclically nonbacktracking walk \emph{bad} if it is anything other
than a repeated walk around a cycle.
%
%
\begin{prop}\label{prop:nooverlaps}
Fix an integer $K$. There is a random time $T$, almost surely finite,
such that there are no bad cyclically nonbacktracking walks of
length $K$ or less in $G(t)$ for all $t\geq T$.
\end{prop}
\begin{pf}
We will work with the discrete-time version of our process $(G_n, n\in
\NN)$.
We first define some machinery
introduced in \cite{LP}.
Consider some cyclically nonbacktracking walk of length $k$
on the edge-labeled complete graph $K_n$ of the form
\[
s_0 \mathop{\longrightarrow}^{w_1} s_1
\mathop{\longrightarrow}^{w_2} s_2\mathop{\longrightarrow}^{w_3}\cdots
\mathop{\longrightarrow}^{w_k}
s_k=s_0.
\]
Here, $s_i\in[n]$ and $w=w_1\cdots w_k$ is the word
of the walk (i.e., each $w_i$
is $\pi_j$ or $\pi_j^{-1}$ for some $j$, indicating which
permutation provided the edge for the walk).
We say that $G_n$ contains the walk if the random permutations
$\pi_1,\ldots,\pi_d$ satisfy $w_i(s_{i-1})=s_i$. In other words,
$G_n$ contains a walk if considering both as edge-labeled directed
graphs, the walk is a subgraph of~$G_n$.

If $(s'_i, 0\leq i\leq k)$ is another walk with the same word, we say
that the two walks are of the same \emph{category} if $s_i=s_j\iff
s'_i=s_j'$. In other words, two walks are of the same category
if they are identical up to relabeling vertices.
The probability that $G_n$ contains a walk depends only on its category.
If a walk contains $e$ distinct edges, then
$G_n$ contains the walk with probability at most $1/[n]_e$.

Let $X_k^{(n)}$ be the number of bad walks of length $k$
in $G_n$ that
start at vertex~$n$.
We will first prove that with probability one, $X_k^{(n)}>0$
for only finitely many $n$.
Call a category bad if the walks in the category are bad.
Let $\Tt_{k,d}$ be the number
of bad categories of walks of length $k$. For any particular
bad category whose walks contain $v$ distinct vertices, there
are $[n-1]_{v-1}$ walks of that category whose first vertex is
$n$. Any bad walk contains
more edges than vertices, so
\[
\E X_k^{(n)} \leq\frac{\Tt_{k,d}[n-1]_{v-1}}{[n]_{v+1}}\leq \frac{\Tt_{k,d}}{n(n-k)}.
\]
Since $X_k^{(n)}$ takes values in the nonnegative integers,
$\P[X_k^{(n)}>0]\leq\E X_k^{(n)}$.
By the Borel--Cantelli lemma, $X_k^{(n)}>0$ for only
finitely many values of $n$.

Thus, for any fixed $r+1$, there exists a random time $N$
such that there are no bad walks on $G_n$ of length $r+1$ or less
starting with vertex $n$,
for $n\geq N$. We claim that for $n\geq N$, there are no bad walks
at all
on $G_n$ with length $r$ or less.
Suppose that $G_m$ contains some bad walk of length $k\leq r$, for some
$m\geq N$. As the graph evolves, it is easy to compute that with
probability one, a~new vertex is eventually inserted into an edge
of this walk.
But at the time $n>m\geq N$ when this occurs, $G_n$ will contain
a bad walk of length $r+1$ or less starting with vertex $n$, a~contradiction. Thus, we have proven that $G_n$ eventually contains
no bad walks of length $r$ or less. The equivalent statement
for the continuous-time version of the graph process follows easily
from this.
\end{pf}

\begin{pf*}{Proof of Theorem \ref{mainthm:eigenvalues}}
Let $\mathrm{CNBW}_{k}^{(s)}(t)$ denote the number of cyclically
nonbacktracking
walks of length $k$ in $G(s+t)$. We decompose these into
those that are repeated walks around cycles of length $j$ for some
$j$ dividing~$k$, and the remaining bad walks, which
we denote $B_{k}^{(s)}(t)$, giving us
\[
\mathrm{CNBW}_{k}^{(s)}(t)=\sum
_{j\mid k}2jC_{j}^{(s)}(t)+
B_{k}^{(s)}(t).
\]

Proposition \ref{prop:nooverlaps} implies that
\[
\lim_{s\to\infty}\P \bigl[\mbox{$B_{k}^{(s)}(t)=0$
for all $k\leq K$, $t\geq0$} \bigr]=1.
\]
By Theorem \ref{mainthm:cycles} together with
the continuous mapping theorem and Slutsky's theorem,
as $s$ tends to infinity,
%
%
%
\begin{equation}
\label{eq:wrongbasis} \bigl(\mathrm{CNBW}_{k}^{(s)}(\cdot),
1\leq k\leq K \bigr)\,{\buildrel\mathcal{L} \over\longrightarrow}\, \biggl( \sum
_{j\divides k}2jN_{j}{({\cdot})}, 1\leq k\leq K \biggr).
\end{equation}

Now, we modify the polynomials $\Gamma_k$
to form a new basis $\{f_k, k\in\NN\}$ with the right properties,
which amounts to expressing each $N_{k}{({t})}$ as a linear combination
of terms $\sum_{j\divides l}2jN_{j}{({t})}$. We do this
with the M\"obius inversion formula. Define the polynomial
%
%
%
\begin{equation}
f_k(x) = \frac{1}{2k}\sum_{j\divides k}
\mu \biggl(\frac{k}{j} \biggr) (2d-1)^{j/2}\Gamma_j(x),
\label{eq:fbasis}
\end{equation}
where $\mu$ is the M\"obius function, given by
\begin{eqnarray*}
\mu(n) = \cases{(-1)^a,&\quad $\mbox{if $n$ is the product of $a$ distinct
primes,}$\vspace*{2pt}
\cr
0,&\quad $\mbox{otherwise.}$}
\end{eqnarray*}
The theorem then follows from \eqref{eq:eigencnbw},
\eqref{eq:wrongbasis}, and the continuous mapping
theorem.
\end{pf*}

%
\begin{pf*}{Proof of Theorem \ref{mainthm:chebycov}} We start by
recalling that, for any fixed $d$,
%
\[
2 \tr  T_i\bigl( G(\infty + t) \bigr)  =
(2d-1)^{-i/2} \sum_{k\mid i} 2k N_k(t).
\]

Now, we will prove finite-dimensional convergence to the stated
Ornstein--Uhlenbeck process.
Consider two time points $s \le t$ and two positive integers~$i, k$. We
will first show that, for any $i, k \in\NN$, the pair $(
(2d-1)^{-i/2}(N_i(s)- \E N_i(s)), (2d-1)^{-k/2}(N_k(t) - \E N_k(t)))$ converges
to a Gaussian limit as $d$ tends to infinity. When $s=t$, this
trivially follows via the central limit theorem and their independent
Poisson joint distribution.

When $s < t$, observe from \eqref{eq:decompN} that
\[
N_k(t)= \sum_{j=1}^k
\alpha(j,k) N_j(s) + Z.
\]
Here $\alpha(j,k) N_j(s)$, $j\in[k]$, and $Z$ are independent Poisson
random variables of various means. Moreover, $Z$ is independent of the
history of the process till time~$s$. Under the stationary law, the
vector $(N_j(s), j\in\NN)$ are independent Poisson random variables.
Thus, if $i >k$, then $N_i(s)$ is independent of $N_k(t)$.
Otherwise, by the thinning property of Poisson, $\alpha(i,k) N_i(s)$ is
independent of $(1- \alpha(i,k))N_i(s)$. Therefore, $N_k(t)-\alpha
(i,k)N_i(s)$, $\alpha(i,k)N_i(s)$, and $(1-\alpha(i,k))N_i(s)$ are
three independent Poisson random variable.

By the 
normal approximation to Poisson, we get the appropriate distributional
convergence to corresponding independent Gaussian random variables.
This shows the joint convergence of $(N_i(s), N_k(t))$ to Gaussian
after centering and scaling.

A similar Gram--Schmidt orthogonalization can be carried out for the
case of time points $t_1\le t_2 \le\cdots\le t_m$ and corresponding
positive integers $j_1, j_2, \ldots, j_m$. This proves the joint
Gaussian convergence of any finite collection of $(N_{j_i}(t_i), i
\in[m])$ under centering and suitable scaling. Since the traces of
Chebyshev polynomials are linear combinations of coordinates of $N$,
the joint Gaussian convergence extends to them by an argument invoking
the continuous mapping and Slutsky's theorems.\vadjust{\goodbreak}

For a fixed $d$, the covariance computation follows from Corollary \ref
{cor:cyclecov} and \eqref{eq:eigencnbw}. Hence, if $s < t$, then
%
%
\begin{eqnarray}
\label{eq:covt}&& \cov\bigl( \tr T_i \bigl( G(\infty+ t) \bigr), \tr
T_j \bigl( G(\infty + s) \bigr) \bigr)
\nonumber
\\[-8pt]
\\[-8pt]
\nonumber
&&\qquad= \frac{1}{4} ( 2d-1
)^{-(i+j)/2}\sum_{k\mid i, l \mid j} 4 lk \cov\bigl(
\limitN_k(t), \limitN_l(s) \bigr).
\end{eqnarray}
Here
\[
\cov \bigl(N_k(t), N_l(s) \bigr) = \cases{
\displaystyle\frac{a(d,l)}{2l}\pmatrix{k-1
\cr
k-l} p^{l} (1-p)^{k-l}
,&\quad $\mbox{$p=e^{s-t}$, if $k \ge l$},$\vspace*{2pt}
\cr
0,&\quad $\mbox{otherwise}.$}
\]

We now fix any $i, j, t, s$ and take $d$ to infinity. Any term $a(d,r)$
is asymptotically the same as $(2d-1)^r$. Thus, the highest order term
(in $d$) on the right-hand side of \eqref{eq:covt} is $(2d-1)^{\min(i,
j)}$. Unless $i=j$, this term is negligible compared to
$(2d-1)^{(i+j)/2}$. This shows that the limiting covariance is zero
unless $i=j$.

On the other hand, when $i=j$, every term on the right-hand side of
\eqref{eq:covt} vanishes, except when $k=i= l=j$. Hence,
\[
\lim_{d\rightarrow\infty}\cov \bigl( \tr T_i \bigl( G(\infty+
t) \bigr), \tr T_i \bigl( G(\infty+ s) \bigr) \bigr)=
\frac{1}{4}2i p^{i}= \frac{i}{2}e^{i(s-t)}.
\]

Finally, we prove the process convergence. One simply needs to argue tightness.
Fix a $K\in\NN$ and, for every $d$, consider the process
\[
\bigl( X_k(t):=(2d-1)^{-k/2} \bigl( 2k N_k(t) -
a(d,k) \bigr), k\in[K], t \ge0 \bigr).
\]
We claim that it suffices to show tightness for this process. This
follows, since then, due to unequal scaling, the difference between
this process and the centered and scaled traces go to zero in
probability as $d$ tends to infinity.

 We sketch a proof of tightness for this process; more details appear in
 \cite{naujasRE}.
 Fix $k$ and $d$. Let $Y(t)$ and $Z(t)$ be counting processes starting at $0$.
 Define $Y(t)$ and $Z(t)$ to jump at points of increase and decrease,
 respectively, of $N_k(t)$.
 We then have $N_k(t)=Y(t)-Z(t)$, and it is not difficult to show
 that $Y(t)$ and $Z(t)$
 are both Poisson processes with rate $a(d,k)/2$.
 Scaled by $(2d-1)^{-k/2}$ and normalized, each converges in law
 to Brownian motion.
 Thus for each $k$, we can write $X_k(t)$
 as a sum of processes converging in law to a limit in $C[0,\infty)$,
 and from here one can obtain the desired tightness of
 $(X_k(t),k\in[K],\,t\geq 0)$ in $D_{\mathbb{R}^K}[0,\infty)$.
\end{pf*}

\begin{appendix}\label{app}
\section*{Appendix: A broad Poisson approximation result}
This Appendix provides the proofs of Theorem \ref{thmm:processapprox}
and Corollary \ref{cor:marginals}.
A~less general version of Theorem \ref{thmm:processapprox} can be
found in
\cite{DJPP}, Theorem 11; we show
in Corollary \ref{cor:poiapprox}(i) how it
follows from Theorem \ref{thmm:processapprox}. Our theorem here also improves
the total variation bound from $O((2d-1)^{2r}/n)$
to $O((2d-1)^{2r-1}/n)$. We conjecture that
Theorem \ref{thmm:processapprox} is sharp.

As in the proof of Theorem 11 in \cite{DJPP}, the main tool is
the Stein--Chen method for Poisson
approximation by size-biased couplings as described in \cite{BHJ},
which uses the following idea:
Recall the definition of $(I_\beta, \beta\in\Ii)$ from
Theorem \ref
{thmm:processapprox}.
For each $\alpha\in\Ii$, let $(J_{\beta\alpha}, \beta\in\Ii)$
be distributed
as $(I_\beta, \beta\in\Ii)$ conditioned\vadjust{\goodbreak} on $I_\alpha=1$.
The goal is to construct a coupling of $(I_\beta, \beta\in\Ii)$
and
$(J_{\beta\alpha}, \beta\in\Ii)$ so that the two random vectors
are ``close together''. We hope that for each $\alpha\in\Ii$,
the cycles in $\Ii\setminus\{\alpha\}$ can be partitioned
into two sets $\Ii_\alpha^{-}$ and $\Ii_\alpha^+$ such that
%
%
%
\begin{eqnarray}
J_{\beta\alpha}&\leq& I_\beta,\qquad \mbox{if $\beta\in\Ii_\alpha^-$,}
\label{eq:monotone-}
\\
J_{\beta\alpha}&\geq& I_\beta,\qquad \mbox{if $\beta\in\Ii_\alpha^+$.}
\label{eq:monotone+}
\end{eqnarray}
If this is the case, then one can approximate $(I_\beta, \beta\in
\Ii)$
by a
Poisson process
by calculating $\cov(I_\alpha, I_\beta)$ for every $\alpha,\beta
\in\Ii$,
according to the following proposition.
%
%
\begin{prop}[(Corollary 10.B.1 in \cite{BHJ})]\label{prop:BHJprop}
Suppose that $\I=(I_\alpha, \alpha\in\Ii)$ is a vector
of 0--1 random variables with $\E I_\alpha=p_\alpha$.
Suppose that $(J_{\beta\alpha}, \beta\in\Ii)$ is distributed
as described above, and that for each $\alpha$ there exists
a partition and a coupling of $(J_{\beta\alpha}, \beta\in\Ii)$
with $(I_\beta, \beta\in\Ii)$ such that \eqref{eq:monotone-}
and \eqref{eq:monotone+}
are satisfied.

Let $\Y=(Y_\alpha, \alpha\in\Ii)$ be a vector of independent
Poisson random variables with $\E Y_\alpha=p_\alpha$. Then
%
%
%
\begin{equation}\quad
d_{\mathrm{TV}}(\I, \Y)\leq\sum_{\alpha\in\Ii}p_\alpha^2
+\sum_{\alpha\in\Ii}\sum_{\smash{\beta\in\Ii_\alpha^-}}
\bigl\vert\cov(I_\alpha,I_\beta)\bigr\vert +\sum
_{\alpha\in\Ii}\sum_{\smash{\beta\in\Ii_\alpha^+}}\cov
(I_\alpha,I_\beta). \label{eq:BHJ}
\end{equation}
\end{prop}
%
We introduce two lemmas, whose proofs we will defer to the end
of the Appendix. The first will let us
approximate $\I$ by $\Z$ rather than by $\Y$, and
the second provides a technical bound that we need.
%
%
\begin{lemma}\label{lem:YZ}
Let $\Y=(Y_\alpha, \alpha\in\Ii)$ and $\Z=(Z_\alpha, \alpha
\in\Ii)$
be vectors of independent Poisson random variables.
Then
\[
d_{\mathrm{TV}}(\Y, \Z)\leq\sum_{\alpha\in\Ii}|\E
Y_\alpha-\E Z_\alpha|.
\]
\end{lemma}
%
%
\begin{lemma}\label{lem:ipbounds}
Let $a$ and $b$ be $d$-dimensional vectors with nonnegative integer
components,
and let $\langle a,b \rangle$ denote the standard Euclidean inner product.
\[
\prod_{i=1}^d\frac{1}{[n]_{a_i+b_i}} - \prod
_{i=1}^d\frac{1}{[n]_{a_i}[n]_{b_i}} \leq
\frac{\langle a,b \rangle}{n}\prod_{i=1}^d
\frac{1}{[n]_{a_i+b_i}}.
\]
\end{lemma}

\begin{pf*}{Proof of Theorem \ref{thmm:processapprox}}
We will give the proof in three sections: First, we make the coupling
and show that it satisfies \eqref{eq:monotone-} and \eqref{eq:monotone+}.
Next, we apply Proposition \ref{prop:BHJprop} to approximate $\I$
by $\Y$, a vector of independent Poissons with $\E Y_\alpha=\E
I_\alpha$.
Last, we approximate $\Y$ by $\Z$ to prove the theorem.

If $d>n^{1/2}$ or $r>n^{1/10}$, then $c(2d-1)^{2r-1}/n>1$
for a sufficiently large choice of $c$, and the theorem holds
trivially. Thus, we will assume throughout that $d\leq n^{1/2}$
and $r\leq n^{1/10}$ (the choice of $1/10$ here is completely arbitrary).
The expression $O(f(d, r, n))$ should be interpreted as a function
of $d$, $r$, and $n$
whose absolute value is bounded by $Cf(d,r,n)$ for some absolute constant
$C$, for all $d$, $r$, and $n$ satisfying $2\leq d\leq n^{1/2}$
and $r\leq n^{1/10}$.

\begin{stepo}[(Constructing the coupling)]

Fix some $\alpha\in\Ii$. We will construct a random vector
$(J_{\beta\alpha}, \beta\in\Ii)$ distributed as
$(I_\beta, \beta\in\Ii)$ conditioned on $I_\alpha=1$.
We do this by constructing a random graph $G_n'$
distributed as $G_n$ conditioned to contain the cycle $\alpha$.
Once this is done, we will define
$J_{\beta\alpha}=1\{\mbox{$G_n'$ contains cycle $\beta$}\}$.

Let $\pi_1,\ldots,\pi_d$ be the random permutations that give rise
to $G_n$.
We will alter them to form permutations $\pi_1',\ldots,\pi_d'$, and
we will
construct $G_n'$ from these. Let us first consider what distributions
$\pi_1',\ldots,\pi_d'$ should have. For example, suppose that
$\alpha$ is the cycle
\[
1 \mathop{\longrightarrow}^{\pi_3} 2
\mathop{\longleftarrow}^{\pi_1} 3\mathop{\longrightarrow}^{\pi_3}4
\mathop{\longrightarrow}^{\pi_1}
=1.
\]
Then $\pi_1'$ should be distributed as a uniform random $n$-permutation
conditioned to make $\pi_1'(3)=2$ and $\pi_1'(4)=1$, and
$\pi_3'$ should be distributed as a uniform random $n$-permutation
conditioned to make $\pi_3'(1)=2$ and $\pi_3'(3)=4$, while
$\pi_2'$ should just be a uniform random $n$-permutation.
A~random graph constructed from $\pi_1'$, $\pi_2'$, and $\pi_3'$
will be distributed as $G_n$ conditioned to contain $\alpha$.

We now describe the construction of $\pi_1',\ldots,\pi_d'$.
Suppose $\alpha$ is the cycle
%
%
%
\begin{equation}
\label{eq:alpha}
s_0 \hspace*{1pt}\mathop{\rule[3pt]{25pt}{0.1pt}}^{w_1}\hspace*{1pt} s_1
\hspace*{1pt}\mathop{\rule[3pt]{25pt}{0.1pt}}^{w_2} s_2
\hspace*{1pt}\mathop{\rule[3pt]{25pt}{0.1pt}}^{w_3}\cdots
\hspace*{1pt}\mathop{\rule[3pt]{25pt}{0.1pt}}^{w_k}
s_k=s_0
\end{equation}
with each edge directed according to whether $w_i(s_{i-1})=s_{i}$
or $w_i(s_{i})=s_{i-1}$.
Fix some $1\leq l\leq d$, and suppose that the edge-label $\pi_l$
appears $M$ times in the cycle~$\alpha$. Let
$(a_m,b_m)$ for $1\leq m\leq M$ be these directed edges.
We must construct $\pi_l'$ to have the uniform
distribution conditioned on $\pi_l'(a_m)=b_m$
for $1\leq m\leq M$.

We define a sequence of random transpositions by the following
algorithm: Let $\tau_1$ swap $\pi_l(a_1)$ with $b_1$.
Let $\tau_2$ swap $\tau_1\pi_l(a_2)$ with $b_2$, and so on.
We then define $\pi_l'=\tau_M\cdots\tau_1\pi_l$.
This permutation satisfies $\pi_l'(a_m)=b_m$ for $1\leq m\leq M$,
and it is distributed uniformly, subject to the
given constraints, which can be
proven by induction on each swap.
We now define $G_n'$ from the permutations $\pi_1',\ldots,\pi_d'$
in the usual way. It is defined on the same probability space as
$G_n$, and it is distributed as $G_n$ conditioned to contain $\alpha$,
giving us a random vector $(J_{\beta\alpha}, \beta\in\Ii)$
coupled with $(I_\beta, \beta\in\Ii)$.

Now, we will give a partition $\Ii^-\cup\Ii^+=\Ii\setminus\{\alpha
\}$
satisfying \eqref{eq:monotone-} and \eqref{eq:monotone+}.
Suppose that
$G_n$ contains an edge ${s_i}\mathop{\longrightarrow}\limits^{w_{i+1}}{v}$
with $v\neq s_{i+1}$,
or an edge ${v}\mathop{\longrightarrow}\limits^{w_{i+1}}{s_{i+1}}$ with $v\neq s_i$.
The graph $G_n'$ cannot contain this edge, since it contains
$\alpha$. In fact,
edges of this form are the \emph{only} ones found in $G_n$ but not~$G_n'$:

%
\begin{lemma}\label{lem:coupling}
Suppose there is an edge
${i}\mathop{\longrightarrow}\limits^{\pi_l}{j}$
contained in $G_n$ but not in $G_n'$. Then $\alpha$
contains either an edge ${i}\mathop{\longrightarrow}\limits^{\pi_l}{v}$
with $v\neq j$,
or $\alpha$ contains an edge ${v}\mathop{\longrightarrow}\limits^{\pi_l}{j}$ with
$v\neq i$.
\end{lemma}
\begin{pf}
Suppose
$\pi_l(i)=j$, but $\pi_l'(i)\neq j$. Then $j$ must
have been swapped
when making $\pi'_l$, which can happen only if
$\pi_l(a_m)=j$ or $b_m=j$ for some $m$. In the first case,
$a_m=i$ and
$\alpha$ contains the edge ${i}\mathop{\longrightarrow}\limits^{\pi_l}{b_m}$ with
$b_m\neq j$,
and in the second $\alpha$ contains the edge
${a_m}\mathop{\longrightarrow}\limits^{\pi_l}{j}$ with $a_m\neq i$.
\end{pf}
Define
$\Ii_\alpha^-$ as all cycles in $\Ii$ that contain an edge
${s_i}\mathop{\longrightarrow}\limits^{w_{i+1}}{v}$ with $v\neq s_{i+1}$
or an edge ${v}\mathop{\longrightarrow}\limits^{w_{i+1}}{s_{i+1}}$
with $v\neq s_i$, and define $\Ii_\alpha^+$
to be the rest of $\Ii\setminus\{\alpha\}$.
Since $G_n'$ cannot contain any cycle in $\Ii_\alpha^{-}$,
we have $J_{\beta\alpha}=0$ for all $\beta\in\Ii_\alpha^-$,
satisfying \eqref{eq:monotone-}.
For any $\beta\in\Ii_\alpha^+$, Lemma \ref{lem:coupling}
shows that if $\beta$ appears in $G_n$, it must also
appear in $G_n'$. Hence $J_{\beta\alpha}\geq I_\beta$, and
\eqref{eq:monotone+} is satisfied.
\end{stepo}

\begin{stepo}[(Approximation of $\I$ by $\Y$)]

The conditions of Proposition \ref{prop:BHJprop} are satisfied, and
we need only bound the sums in \eqref{eq:BHJ}.
Let $p_\alpha=\E I_\alpha$, the probability that cycle $\alpha$
appears in $G_n$. Recall that this equals
$\prod_{i=1}^d 1/[n]_{e_i}$, where
$e_i$ is the number of times $\pi_i$ and $\pi_i^{-1}$ appear in the word
of $\alpha$. This means that
%
%
%
\begin{equation}
\frac{1}{n^k}\leq p_\alpha\leq\frac{1}{[n]_k},\label{eq:pbound}
\end{equation}
where $k=\vert\alpha\vert$, the length of cycle $\alpha$.

We bound the first sum in \eqref{eq:BHJ} by
%
%
%
\begin{eqnarray}\label{eq:psquaredcov}
\sum_{\alpha\in\Ii}p_\alpha^2& =&\sum
_{k=1}^r\sum_{\alpha\in\Ii_k}p_\alpha^2
\leq\sum_{k=1}^r\sum
_{\alpha\in\Ii_k}\frac{1}{[n]_k^2}= \sum_{k=1}^r \biggl(
\frac{[n]_ka(d,k)}{2k} \biggr) \biggl( \frac{1}{[n]_k^2} \biggr)
\nonumber
\\[-8pt]
\\[-8pt]
\nonumber
&\leq&\sum_{k=1}^r \frac{2d(2d-1)^{k-1}}{2k[n]_k}= O
\biggl(\frac{d}{n} \biggr).\nonumber 
\end{eqnarray}

To bound the second sum in \eqref{eq:BHJ},
we investigate the size of $\Ii_\alpha^-$.
Suppose that $\alpha\in\Ii_k$, and $\alpha$ has the form
given in \eqref{eq:alpha}. Any $\beta\in\Ii_\alpha^-$ must
contain an edge
${s_i}\mathop{\longrightarrow}\limits^{w_{i+1}}{v}$ with $v\neq s_{i+1}$,
or an edge ${v}\mathop{\longrightarrow}\limits^{w_{i+1}}{s_{i+1}}$ with $v\neq s_{i}$, and
there are at most $2k(n-1)$ edges of this form.
For any given edge, there are at most $[n-2]_{j-2}(2d-1)^{j-1}$
cycles in $\Ii_j$ that contain that edge, for any $j\geq2$.
Thus for any $\alpha\in\Ii_k$, the number of cycles of
length $j\geq2$ in $\Ii_\alpha^-$ is at most
$2k[n-1]_{j-1}(2d-1)^{j-1}$, and this bound
also holds for $j=1$.

For any $\beta\in\Ii_\alpha^-$, it holds that $\E[I_\alpha I_\beta]=0$,
so that $\cov(I_\alpha,I_\beta)=-p_\alpha p_\beta$. Putting this all
together and applying \eqref{eq:pbound}, we have
%
%
%
\begin{eqnarray}\label{eq:-cov}
\sum_{\alpha\in\Ii}\sum_{\smash{\beta\in\Ii_\alpha^-}}
\bigl\vert\cov(I_\alpha,I_\beta)\bigr\vert &= &\sum
_{k=1}^r\sum_{\alpha\in\Ii_k}\sum
_{j=1}^r\sum_{\beta\in
\Ii
_\alpha^-\cap\Ii_j}
p_\alpha p_\beta
\nonumber
\\
&\leq&\sum_{k=1}^r\vert\Ii_k
\vert \frac{1}{[n]_k} \sum_{j=1}^r\bigl\vert
\Ii_\alpha^-\cap\Ii_j\bigr\vert \frac{1}{[n]_j}
\\
&\leq&\sum_{k=1}^r\frac{a(d,k)}{2k}\sum
_{j=1}^r\frac{2k(2d-1)^{j-1}}{n}
\nonumber
\\
&=&\sum_{k=1}^ra(d,k)O \biggl(
\frac{(2d-1)^{r-1}}{n} \biggr) =O \biggl(\frac{(2d-1)^{2r-1}}{n}
\biggr).\nonumber
\end{eqnarray}

The final sum in \eqref{eq:BHJ} is the most difficult
to bound. We partition
$\Ii_\alpha^+$ into sets
$\Ii_\alpha^+=\Ii_\alpha^0\cup\cdots\cup\Ii_\alpha^{\vert
\alpha\vert-1}$,
where $\Ii_\alpha^l$ is all cycles in $\Ii_\alpha^+$
that share exactly $l$ labeled edges with $\alpha$.
For any $\beta\in\Ii_\alpha^+$,
\[
\E[I_\alpha I_\beta] = \P[\mbox{$G$ contains $\alpha$ and $
\beta$}] =\prod_{i=1}^d\frac{1}{[n]_{e_i}},
\]
where $e_i$ is the number of $\pi_i$-labeled edges
in $\alpha\cup\beta$.
Thus for $\beta\in\Ii_\alpha^l$,
%
%
%
\begin{equation}
\frac{1}{n^{\vert\alpha\vert+\vert\beta\vert-l}}\leq\E [I_\alpha I_\beta]\leq
\frac{1}{[n]_{\vert\alpha\vert+\vert\beta\vert-l}}. \label{eq:covbound}
\end{equation}

We start by seeking estimates on the size
of~$\Ii_\alpha^l$ for $l\geq1$.
Fix some choice of $l$ edges
of $\alpha$. We start by counting the cycles in $\Ii_\alpha^l$
that share exactly these edges with $\alpha$.
We illustrate this in Figure~\ref{fig:graphassembly}.
Call the graph consisting of these edges
$H$, and suppose that $H$ has $p$ components.
Since it is a forest, $H$ has $l+p$ vertices.
\begin{figure}
\centering
\begin{tabular}{@{}cc@{}}

\includegraphics{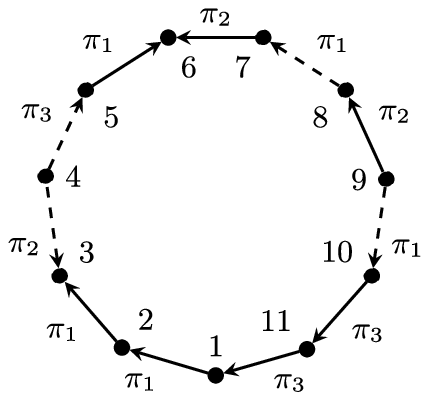}
 & \includegraphics{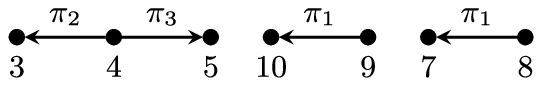}\\
\footnotesize{The cycle $\alpha$, with $H$ dashed.} &
\footnotesize{\emph{Step} 1. We lay out the components}\\
\footnotesize{The subgraph
$H$ has components}&\footnotesize{$A_1,\ldots,A_p$.
We can order and orient}\\
\footnotesize{$A_1,\ldots,A_p$. In this example, the}&\footnotesize{$A_2,\ldots,A_p$ however we would like,
for}\\
\footnotesize{number of components of $H$ is}&\footnotesize{a total of $(p-1)!2^{p-1}$ choices.
Here, we}\\
\footnotesize{$p=3$,
the size of $\alpha$ is $k=11$, and}&\footnotesize{have ordered the components $A_1, A_3, A_2$,}\\
\footnotesize{the number
of edges in $H$ is $l=4$.}&\footnotesize{and we have reversed the orientation}\\
\footnotesize{In this example, we will construct}&\footnotesize{of $A_3$.}\\
\footnotesize{a cycle $\beta$ of length
$j=10$
that}&\footnotesize{}\\
\footnotesize{overlaps with $\alpha$ at $H$.}&\footnotesize{}\\[3pt]

\includegraphics{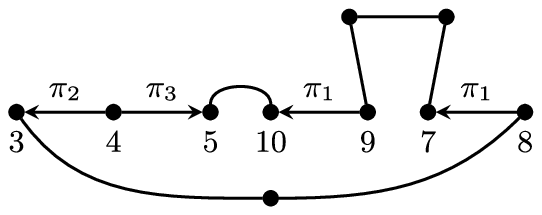}
 & \includegraphics{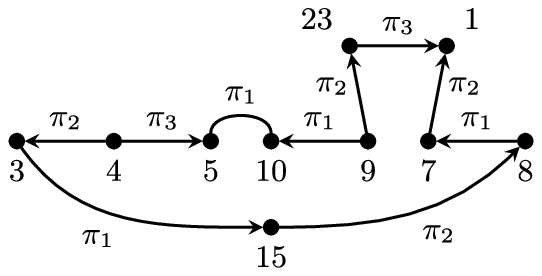}\\
\footnotesize{\emph{Step} 2. Next, we choose how many edges} &
\footnotesize{\textit{Step} 3. We can choose the new vertices in}\\
\footnotesize{will go in each
gap between components.}&\footnotesize{$[n-p-l]_{j-p-l}$ ways, and we can direct}\\
\footnotesize{Each gap must contain at least
one edge,}&\footnotesize{and give labels
to the new edges in at}\\
\footnotesize{and we must add a total of $j-l$ edges,}&\footnotesize{most $(2d-1)^{j-l}$ ways.}\\
\footnotesize{giving
us ${j-l-1\choose p-1}$ choices.}&\footnotesize{}\\
\footnotesize{In this example, we have added one edge}&\footnotesize{}\\
\footnotesize{after $A_1$,
three after $A_3$, and two after $A_2$.}&\footnotesize{}\\
\end{tabular}
\caption{Assembling an element $\beta\in\Ii_\alpha^l$ that
overlaps with $\alpha$ at a given subgraph $H$.}
\label{fig:graphassembly}
\end{figure}

Let $A_1,\ldots, A_p$ be the components of $H$. We can
assemble any element $\beta\in\Ii_\alpha^l$
that overlaps with $\alpha$ in $H$
by stringing together these components in some order, with
other edges in between.
Each component can appear in
$\beta$ in one of two orientations. Since the vertices in $\beta$
have no fixed ordering, we can
assume without
loss of generality that
$\beta$ begins with component $A_1$ with a fixed orientation.
This leaves $(p-1)!2^{p-1}$ choices for the order
and orientation of $A_2,\ldots,A_p$ in $\beta$.

Imagine now the components laid out in a line,
with gaps between them, and count the number
of ways to fill the gaps. Suppose that $\beta$
is to have length $j$.
Each of the
$p$ gaps
must contain at least one edge, and the total number
of edges in all the gaps is $j-l$.
Thus, the total number of possible gap sizes
is the number of compositions of $j-l$ into $p$ parts,
or ${j-l-1\choose p-1}$.

Now that we have chosen the number of edges
to appear in each gap, we choose the edges themselves.
We can do this by giving
an ordered list $j-p-l$ vertices to go in the gaps,
along with a label and an orientation
for each of the $j-l$ edges this gives.
There are $[n-p-l]_{j-p-l}$ ways to choose the vertices.
We can give each new edge any orientation and label subject
to the constraint that the word of the cycle we
construct must be reduced.
This means we have at most $2d-1$ choices for the orientation
and label of each new edge, for a total of at most $(2d-1)^{j-i}$.

All together, there are at most
$(p-1)!2^{p-1}{j-l-1\choose p-1}[n-p-l]_{j-p-l}(2d-1)^{j-l}$
elements of $\Ii_j$ that overlap with the cycle $\alpha$ at the
subgraph $H$. We now calculate the number of different ways to choose
a subgraph $H$ of $\alpha$ with $l$ edges and $p$ components.
Suppose $\alpha$ is given as in \eqref{eq:alpha}.
We first choose a vertex $s_{i_0}$. Then, we can
specify which edges to include in $H$ by giving a sequence
$a_1,b_1,\ldots,a_p,b_p$ instructing us
to include in $H$ the first $a_1$ edges after $s_{i_0}$,
then to exclude the next~$b_1$, then to include the next $a_2$,
and so on. Any sequence for which $a_i$ and $b_i$
are positive integers, $a_1+\cdots+ a_p=l$,
and $b_1+\cdots+b_p=k-l$ gives us a valid choice of $l$ edges of
$\alpha$
making up $p$ components. This counts each subgraph $H$ a total of
$p$ times, since we could begin with any component of $H$.
Hence, the number of subgraphs $H$ with $l$ edges
and $p$ components is $(k/p){l-1\choose p-1}{k-l-1\choose p-1}$.
This gives us the bound
\begin{eqnarray*}
\bigl|\Ii_\alpha^l\cap\Ii_j\bigr|&\leq&\sum
_{p=1}^{l\wedge(j-l)} (k/p)\pmatrix{l-1
\cr
p-1}\pmatrix{k-l-1
\cr
p-1} (p-1)!
\\
&&\hspace*{27pt}{}\times 2^{p-1}\pmatrix{j-l-1
\cr
p-1}[n-p-l]_{j-p-l}(2d-1)^{j-l}.
\end{eqnarray*}
We apply the bounds
\begin{eqnarray*}
\pmatrix{l-1
\cr
p-1}&\leq&\frac{r^{p-1}}{(p-1)!},
\\
\pmatrix{k-l-1
\cr
p-1}, \pmatrix{j-l-1
\cr
p-1}&\leq& \bigl(er/(p-1)
\bigr)^{p-1},
\end{eqnarray*}
to get
\begin{eqnarray*}
\bigl|\Ii_\alpha^l\cap\Ii_j\bigr| &\leq&
k(2d-1)^{j-l}[n-1-l]_{j-1-l}  \\
&&{}\times\Biggl(1+ \sum
_{p=2}^{i\wedge(k-i)}\frac{1}{p} \biggl(
\frac{2e^2r^3}{(p-1)^2} \biggr)^{p-1} \frac{1}{[n-1-l]_{p-1}} \Biggr).
\end{eqnarray*}
Since $r\leq n^{1/10}$,
the sum in the above equation is bounded by an absolute constant.
Applying this bound and \eqref{eq:covbound},
for any $\alpha\in\Ii_k$ and $l\geq1$,
\begin{eqnarray*}
\sum_{\beta\in\Ii_\alpha^l} \cov(I_\alpha,I_\beta)
&\leq&\sum_{j=l+1\vphantom{\Ii_\alpha^l}}^r \sum
_{\beta\in\Ii_\alpha^l\cap\Ii_j} \frac{1}{[n]_{k+j-l}}
\leq\sum_{j=l+1}^rO \biggl(
\frac{k(2d-1)^{j-l}}{n^{k+1}} \biggr)
\\
&=&O \biggl(\frac{k(2d-1)^{r-l}}{n^{k+1}} \biggr).
\end{eqnarray*}
Therefore,
%
%
%
\begin{eqnarray}\label{eq:+cov}
\sum_{\alpha\in\Ii}\sum_{l\geq1}
\sum_{\beta\in\Ii_\alpha^l}\cov(I_\alpha,I_\beta)
&=& \sum_{k=1}^r\sum
_{\alpha\in\Ii_k}\sum_{l=1}^{k-1}
\sum_{\beta\in\Ii_\alpha^l}\cov(I_\alpha,I_\beta)
\nonumber
\\
&\leq&\sum_{k=1}^r\sum
_{\alpha\in\Ii_k}\sum_{l=1}^{k-1} O
\biggl(\frac{k(2d-1)^{r-l}}{n^{k+1}} \biggr)
\\
&=& \sum_{k=1}^r\frac{[n]_ka(d,k)}{2k} O
\biggl(\frac{k(2d-1)^{r-1}}{n^{k+1}} \biggr)
\nonumber\\
&= &\sum_{k=1}^r O \biggl(
\frac{(2d-1)^{r+k-1}}{n} \biggr)
\nonumber
\\
&= &O \biggl(\frac{(2d-1)^{2r-1}}{n} \biggr).\nonumber
\end{eqnarray}

Last, we must bound $\sum_{\alpha\in\Ii}\sum_{\beta\in\Ii
_\alpha^0}\cov
(I_\alpha,
I_\beta)$.
For any word $w$, let
$e^w_i$ be the number
of appearances of $\pi_i$ and $\pi_i^{-1}$ in $w$.
Let $\alpha$ and $\beta$ be cycles with words $w$
and $u$, respectively, and let $k=\vert\alpha\vert$ and $j=\vert
\beta\vert$.
Suppose that $\beta\in\Ii^0_\alpha$. Then
\begin{eqnarray*}
\cov(I_\alpha,I_\beta)&=& \prod_{i=1}^d
\frac{1}{[n]_{e^w_i+e^u_i}} - \prod_{i=1}^d
\frac{1}{[n]_{e^w_i}[n]_{e^u_i}}
\\
&\leq&\frac{\langle e^w,e^u \rangle}{n}\prod_{i=1}^d
\frac{1}{[n]_{e^w_i+e^u_i}} \leq\frac{\langle e^w,e^u \rangle}{n[n]_{k+j}}
\end{eqnarray*}
by Lemma \ref{lem:ipbounds}.
For any pair of words $w\in\Ww_k$ and $u\in\Ww_j$, there are
at most
$[n]_{k}[n]_j$ pairs of cycles $\alpha,\beta\in\Ii$ with words
$w$ and $u$, respectively. Enumerating over all
$w\in\Ww_k$ and $u\in\Ww_j$, we count each pair
of cycles $\alpha,\beta$ exactly $4kj$ times. Thus,
\begin{eqnarray*}
\sum_{\alpha\in\Ii_k}\sum_{\beta\in\Ii_\alpha^0\cap\Ii
_j}
\cov(I_\alpha,I_\beta) &\leq&\frac{[n]_{k}[n]_j}{4kjn[n]_{k+j}}\sum
_{w\in\Ww_k}\sum_{u\in
\Ww
_j} \bigl\langle
e^w,e^u \bigr\rangle
\\
&\leq&\frac{1+O(r^2/n)} {
4kjn}\biggl\langle\sum_{w\in\Ww_k}e^w,
\sum_{u\in\Ww_j}e^u \biggr\rangle.
\end{eqnarray*}
The vector $\sum_{w\in\Ww_k}e^w$ has every entry equal by symmetry,
as does $\sum_{u\in\Ww_j}e^u$. Thus, each entry of
$\sum_{w\in\Ww_k}e^w$ is $ka(d,k)/d$, and each entry
of $\sum_{u\in\Ww_j}e^u$
is $ja(d,j)/d$. The inner product in the above equation
comes to\break  $kja(d,k)a(d,j)/d$, giving us
\begin{eqnarray*}
\sum_{\alpha\in\Ii_k}\sum_{\beta\in\Ii_\alpha^0\cap\Ii
_j}
\cov(I_\alpha,I_\beta) &\leq&\frac{a(d,k)a(d,j)(1+O(r^2/n))}{4dn}
\\
&=& O \biggl(\frac{(2d-1)^{j+k-1}}{n} \biggr).
\end{eqnarray*}
Summing over all $1\leq k,j\leq r$,
%
%
%
\begin{equation}
\sum_{\alpha\in\Ii}\sum_{\beta\in\Ii^0_\alpha}
\cov(I_\alpha,I_\beta) = \biggl(\frac{(2d-1)^{2r-1}}{n} \biggr).
\label{eq:0cov}
\end{equation}

We can now combine equations \eqref{eq:psquaredcov}, \eqref{eq:-cov},
\eqref{eq:+cov}, and \eqref{eq:0cov} with
Proposition \ref{prop:BHJprop} to show that
%
%
%
\begin{equation}
d_{\mathrm{TV}}(\I, \Y)=O \biggl(\frac{(2d-1)^{2r-1}}{n} \biggr).\label{eq:IY}
\end{equation}
\end{stepo}

\begin{stepo}[(Approximation of $\Y$ by $\Z$)]

By Lemma \ref{lem:YZ} and \eqref{eq:pbound},
\begin{eqnarray*}
d_{\mathrm{TV}}(\Y, \Z)\leq\sum_{\alpha\in\Ii}\vert\E
Y_\alpha-\E Z_\alpha\vert &\leq&\sum
_{k=1}^r \sum_{\alpha\in\Ii_k}
\biggl( \frac{1}{[n]_k}- \frac{1}{n^k} \biggr)
\\
&=&\sum_{k=1}^r\frac{a(d,k)}{2k}
\biggl(1-\frac{[n]_k}{n^k} \biggr).
\end{eqnarray*}
Since $[n]_k\geq n^k(1-k^2/2n)$,
\[
d_{\mathrm{TV}}(\Y, \Z)\leq\sum_{k=1}^r
\frac{a(d,k)k}{4n} =O \biggl(\frac{r(2d-1)^r}{n} \biggr).
\]
Together with \eqref{eq:IY}, this bounds the total variation
distance between the laws of $\I$ and $\Z$ and proves the
theorem.\qed
\end{stepo}
\noqed\end{pf*}
The distributions of any functionals of $\I$ and $\Z$ satisfy the
same bound in total variation distance. This gives us several
results as easy corollaries, including an improvement
on \cite{DJPP}, Theorem 11.
%
%
\begin{cor}\label{cor:poiapprox}
\textup{(i)} 
Let $(Z_k, 1\leq k\leq r)$ be a vector
of independent Poisson random variables with $\E Z_k=a(d,k)/2k$.
Let $C_k$ denote the number of $k$-cycles in $G_n$, a $2d$-regular
permutation random graph on $n$ vertices. Then for some absolute
constant $c$,
\[
d_{\mathrm{TV}} \bigl((C_k, 1\leq k\leq r), (Z_k, 1
\leq k\leq r) \bigr) \leq\frac{c(2d-1)^{2r-1}}{n}.
\]

\begin{longlist}[(ii)]
\item[(ii)] Let $(Z_w, w\in\Ww'_K)$ be a vector \textup{(ii)}
of independent Poisson random variables with $\E Z_w=1/h(w)$.
Let $C_w$ denote the number of cycles with word $w$
in $G_n$, a $2d$-regular
permutation random graph on $n$ vertices. Then for some absolute
constant $c$,
\[
d_{\mathrm{TV}} \bigl( \bigl(C_w, w\in\Ww'_K
\bigr), \bigl(Z_w, w\in\Ww'_K \bigr)
\bigr) \leq \frac{c(2d-1)^{2K-1}}{n}.
\]
\end{longlist}
\end{cor}
\begin{pf}
Observe that $C_k=\sum_{\alpha\in\Ii_k} I_\alpha$,
and that if we define $Z_k=\break \sum_{\alpha\in\Ii_k} Z_\alpha$,
then $(Z_k, 1\leq k\leq r)$ is distributed as described.
Thus (i) follows from Theorem~\ref{thmm:processapprox}.

To prove \textup{(ii)}, note that $C_w=\sum_\alpha I_\alpha$,
where the sum is over all cycles in $\Ii$ with word $w$.
We then define $Z_w$ as the analogous sum over $Z_\alpha$.
Since the number of cycles in $\Ii$ with word $w$ is $[n]_k/h(w)$,
we have $\E Z_w=1/h(w)$, and the total variation bound
follows from Theorem \ref{thmm:processapprox}.
\end{pf}
We can also use Theorem \ref{thmm:processapprox} to bound the likelihood
that $G_n$ contains two overlapping cycles of size $r$ or less.
%
%
\begin{cor}\label{cor:overlap}
Let $G_n$ be a $2d$-regular permutation random graph on $n$
vertices. Let $E$ be the event that $G_n$ contains two cycles
of length $r$ or less with a vertex in common.
Then for some absolute constant $c'$,
for all $d\geq2$ and $n,r\geq1$,
\[
\P[E]\leq\frac{c'(2d-1)^{2r}}{n}.
\]
\end{cor}
\begin{pf}
Let $E'$ be the event that $Z_\alpha=Z_\beta=1$ for two cycles
$\alpha, \beta\in\Ii$ that have a vertex in common.
By Theorem \ref{thmm:processapprox},
\[
\P[E]\leq\P \bigl[E' \bigr]+\frac{c(2d-1)^{2r-1}}{n}.
\]
For any cycle $\alpha\in\Ii_k$, there are at most
$k[n-1]_{j-1}a(d,j)$ cycles in $\Ii_j$ that share a vertex
with $\alpha$. For any such cycle $\beta$, the chance
that $Z_\alpha=1$ and $Z_\beta=1$ is less than $1/[n]_k[n]_j$.
By a union bound,
\begin{eqnarray*}
\P \bigl[E' \bigr]&\leq&\sum_{k=1}^r
\frac{a(d,k)[n]_k}{2k}\sum_{j=1}^r
\frac{k[n-1]_{j-1}a(d,j)}{[n]_k[n]_j}
\\
&\leq&\sum_{k=1}^r\sum
_{j=1}^r\frac{a(d,k)a(d,j)}{2n}=O \biggl(
\frac{(2d-1)^{2r}}{n} \biggr).
\end{eqnarray*}
\upqed\end{pf}


\begin{pf*}{Proof of Corollary \ref{cor:marginals}}
When $d=1$, there is only one word of each length in $\Ww'_K$,
and statement (i) reduces to the well-known fact that
the cycle
counts of a random permutation converge to independent Poisson random
variables (see \cite{AT} for much more on this subject).
In this case, $G(t)$ is made up of disjoint cycles for all times $t$,
so that statement (ii) is trivially satisfied.

When $d\geq2$, let $C_{w}^{(n)}$ be the number of cycles with word $w$ in
$G_n$,
as in Corollary \ref{cor:poiapprox}(ii). The random vector
$(C_w(t), w\in\Ww'_K)$ is a mixture of the random
vectors $(C_{w}^{(n)}, w\in\Ww'_K)$ over different values of $n$.
That is,
\[
\P \bigl[ \bigl(C_w(t), w\in\Ww'_K \bigr)
\in A \bigr] =\sum_{n=1}^{\infty}
\P[M_t=n]\P \bigl[ \bigl(C_{w}^{(n)}, w\in
\Ww'_K \bigr)\in A \bigr]
\]
for any set $A$, recalling that $G(t)=G_{M_t}$.
Corollary \ref{cor:poiapprox}(ii) together
with the fact that $\P[M_t>N]\to1$ as $t\to\infty$ for
any $N$ imply that $(C_w(t), w\in\Ww'_K)$ converges
in law to $(Z_w, w\in\Ww'_K)$, establishing
statement (i). Statement (ii)
follows in the same way from Corollary \ref{cor:overlap}.
\end{pf*}
\begin{pf*}{Proof of Lemma \ref{lem:YZ}}
We will apply the Stein--Chen method directly.
Define the operator $\Aa$ by
\[
\Aa h(x) = \sum_{\alpha\in\Ii}\E[Z_\alpha]
\bigl(h(x+e_\alpha)-h(x) \bigr) +\sum_{\alpha\in\Ii}x_\alpha
\bigl(h(x-e_\alpha)-h(x) \bigr)
\]
for any $h\dvtx \ZZ_+^{\vert\Ii\vert}\to\RR$ and $x\in\ZZ
_+^{\vert\Ii\vert}$.
This is the Stein operator for the law of $\Z$, and
$\E\Aa h(\Z)=0$ for any bounded function $h$.
By Proposition 10.1.2 and Lemma~10.1.3 in \cite{BHJ},
for any set $A\subseteq\ZZ_+^{\vert\Ii\vert}$, there is a function
$h$ such that
\[
\Aa h(x)=1\{x\in A\}- \P[\Z\in A],
\]
and this function has the property that
%
%
%
\begin{equation}
\mathop{\sup_{x\in\ZZ_+^{\vert\Ii\vert}}}_{ \alpha\in\Ii
}\bigl|h(x+e_\alpha
)-h(x)\bigr|\leq1. \label{eq:mvsteinbound1}
\end{equation}
Thus we can bound the total variation distance between
the laws of $\Y$ and $\Z$ by bounding $\vert\E\Aa h(\Y)\vert$
over all such functions $h$.

We write $\Aa h(\Y)$ as
\begin{eqnarray*}
\Aa h(\Y) &=& \sum_{\alpha\in\Ii}\E[Y_\alpha]
\bigl(h(\Y+e_\alpha)-h(\Y) \bigr) +\sum_{\alpha\in\Ii}Y_\alpha
\bigl(h(\Y-e_\alpha)-h(x) \bigr)
\\
&&{} +\sum_{\alpha\in\Ii} (\E Z_\alpha- \E
Y_\alpha) \bigl(h(\Y+e_\alpha)-h(\Y) \bigr).
\end{eqnarray*}
The first two of these sums have expectation zero, so
\[
\bigl\vert\E\Aa h(\Y)\bigr\vert \leq\sum_{\alpha\in\Ii}\vert\E
Z_\alpha- \E Y_\alpha\vert \E\bigl\vert h(\Y+e_\alpha)-h(
\Y)\bigr\vert.
\]
By \eqref{eq:mvsteinbound1},
$\vert h(\Y+e_\alpha)-h(\Y)\vert\leq1$,
which proves the lemma.
\end{pf*}
%
%
\begin{pf*}{Proof of Lemma \ref{lem:ipbounds}}
We define a family of independent
random maps $\sigma_i$ and $\tau_i$
for $1\leq i\leq d$.
Choose $\sigma_i$ uniformly from all injective maps from $[a_i]$
to~$[n]$, and choose $\tau_i$ uniformly from all injective maps
from $[b_i]$ to $[n]$. Effectively, $\sigma_i$ and $\tau_i$
are random ordered subsets of $[n]$.
We say that $\sigma_i$ and $\tau_i$ \textit{clash} if their images
overlap.
\[
\P[\mbox{$\sigma_i$ and $\tau_i$ clash for some $i$}]
= 1-\prod_{i=1}^d\frac{[n]_{a_i+b_i}}{[n]_{a_i}[n]_{b_i}}.
\]
For any $1\leq i\leq d$, $1\leq j\leq a_i$, and $1\leq k\leq b_i$,
the probability that $\sigma_i(j)=\tau_i(k)$ is $1/n$.
By a union bound,
\[
\P[\mbox{$\sigma_i$ and $\tau_i$ clash for some $i$}]
\leq\sum_{i=1}^d\frac{a_ib_i}{n}=
\frac{\langle a,b \rangle}{n}.
\]
We finish the proof by dividing
both sides of this inequality by $\prod_{i=1}^d[n]_{a_i+b_i}$.
\end{pf*}
\end{appendix}

%
%

%



\printaddresses

\end{document}